\RequirePackage{fix-cm}
\documentclass[smallextended]{svjour3}       
\smartqed  
\usepackage{graphicx}
\RequirePackage{amsmath,amsfonts,amssymb,mathtools,thmtools}
\RequirePackage[round]{natbib}
\RequirePackage[colorlinks,citecolor=blue,urlcolor=blue]{hyperref}
\RequirePackage{enumitem,graphicx,mathrsfs}

\usepackage{relsize}
\usepackage{epstopdf}
\usepackage{bm}
\usepackage{todonotes}
\usepackage{epigraph} 
\usepackage{graphicx,psfrag,epsf} 
\usepackage{subcaption}
\usepackage{nicefrac}

\newcommand{\refFig}[1]{Fig.~{\ref{#1}}}
\newcommand{\refEq}[1]{Eq.~{\eqref{#1}}}

\newcommand{\refThm}[1]{Theorem~{\ref{#1}}}

\newcommand{\refApp}[1]{Appendix~{\ref{#1}}}

\newcommand{\bibOrder}[1]{}

\newcommand{\field}[1]{\mathbb{#1}} 
\newcommand{\R}{\field{R}} 

\newcommand{\Hc}{\mathcal{H}}

\newcommand{\Nc}{\mathcal{N}}

\newcommand{\Oc}{\mathcal{O}}
\newcommand{\Pc}{\mathcal{P}}
\newcommand{\Pf}{\mathbb{P}}

\newcommand{\der}{\mathrm{d}}

\newcommand{\iid}{\textnormal{iid}}

\newcommand{\iidSim}{\overset{\iid}{\sim}}

\newcommand{\cov}{\textnormal{Cov}}

\newcommand{\BF}{\textnormal{BF}}

\newcommand{\revise}[1]{\textcolor{black}{#1}}

\definecolor{nRed}{RGB}{220,0,12}
\definecolor{nOrange}{RGB}{254,166,0}
\definecolor{nBlue}{RGB}{40,96,180}

\begin{document}

\title{Bayes Factors for Peri-Null Hypotheses
}


\author{Alexander Ly         \and
        Eric-Jan Wagenmakers  
}


\institute{A. Ly \at
              University of Amsterdam\\
              Centrum Wiskunde \& Informatica \\
              \email{a.ly@uva.nl}           
           \and
           E.--J. Wagenmakers \at
           University of Amsterdam\\
           \email{ejwagenmakers@gmail.com}           
}

\date{Received: date / Accepted: date}

\maketitle

\begin{abstract}
A perennial objection against Bayes factor point-null hypothesis tests is that the point-null hypothesis is known to be false from the outset. We examine the consequences of approximating the sharp point-null hypothesis by a hazy `peri-null' hypothesis instantiated as a narrow prior distribution centered on the point of interest. The peri-null Bayes factor then equals the point-null Bayes factor multiplied by a correction term which is itself a Bayes factor. For moderate sample sizes, the correction term is relatively inconsequential; however, for large sample sizes the correction term becomes influential and causes the peri-null Bayes factor to be inconsistent and approach a limit that depends on the ratio of prior ordinates evaluated at the maximum likelihood estimate. We characterize the asymptotic behavior of the peri-null Bayes factor and briefly discuss suggestions on how to construct peri-null Bayes factor hypothesis tests that are also consistent.

\keywords{Consistency \and Peri-null correction factor \and Asymptotic sampling distribution}
\subclass{62F03 \and 62F15 \and 62F05}
\end{abstract}	

\newpage

\setlength{\epigraphwidth}{0.4\textwidth} 
\epigraph{vagueness leads nowhere.}{Jeffreys, 1937}
\setlength{\epigraphwidth}{0.4\textwidth} 
\nocite{Jeffreys1937Aristotelian} 

\section{Introduction}
\label{secIntro}
In the Bayesian paradigm, the support that data \( y^{n}:= (y_{1} , \ldots, y_{n}) \) provide for an alternative hypothesis \( \Hc_{1} \) versus a point-null hypothesis \( \Hc_{0} \) is given by the Bayes factor \( \BF_{10}(y^{n}) \):
\begin{align}
\underbrace{\frac{p(y^{n} \mid \mathcal{H}_{1})}{p(y^{n} \mid \mathcal{H}_{0})}}_{\BF_{10} (y^{n})} \quad 
&=  \overbrace{\frac{P(\mathcal{H}_{1} \, | \,  y^{n})}{P(\mathcal{H}_{0} \, | \,  y^{n})}}^{\text{Posterior model odds}} \, \bigg/ \overbrace{\frac{ P(\mathcal{H}_{1})}{P(\mathcal{H}_{0})}}^{\text{prior model odds}} \\
&= \frac{\int_{\Theta_{1}} f(y^{n} \, | \, \theta_{1}) \, \pi(\theta_{1} \, | \, \Hc_{1}) \, \der \theta_{1}}{\int_{\Theta_{0}} f(y^{n} \, | \, \theta_{0}) \, \pi(\theta_{0} \, | \, \Hc_{0}) \, \der \theta_{0}},
\end{align}    
where the first line indicates that the Bayes factor quantifies the change from prior to posterior model odds \citep{WrinchJeffreys1921}, and the second line indicates that this change is given by a ratio of marginal likelihoods, that is, a comparison of prior predictive performance obtained by integrating the parameters \( \theta_{j} \) out of the \( j \)th model's likelihood \( f( y^{n} \, | \, \theta_{j}) \) at the observations \( y^{n} \) with respect to the prior density \( \pi(\theta_{j} \, | \, \Hc_{j}) \) \citep{Jeffreys1935,Jeffreys1939,KassRaftery1995}.  Although the general framework applies to the comparison of any two models (as long as the models make probabilistic predictions; \citealp{Dawid1984, shafer2019game}), the procedure developed by Harold Jeffreys in the late 1930s was explicitly designed as an improvement on \( p \)-value null-hypothesis significance testing. 

In the prototypical scenario, a null hypothesis \( \Hc_{0} \) has \( p_{0} \) free parameters, whereas an alternative hypothesis \( \Hc_{1} \) has \( p = p_{0} + 1 \) free parameters; the additional free parameter in \( \Hc_{1} \) is the one that is \emph{test-relevant}. For instance, in Jeffreys's \( t \)-test the test-relevant parameter \( \delta = \mu / \sigma \) represents standardized effect size; after assigning prior distributions to the models' parameters we may compute the Bayes factor in favor of \( \Hc_{0}: \delta = 0 \) with free parameter \( \theta_{0} =  \sigma \in (0, \infty) \) against \( \Hc_{1}: \theta_{1} =  ( \delta, \sigma ) \in \R \times (0, \infty) \) where \( \delta  \) is unrestricted and where \( \sigma \) denotes the common nuisance parameter. When \( \BF_{01}(y^{n}) = 1/\BF_{10}(y^{n}) \) is larger than 1, the data provide evidence that the `general law' \( \Hc_{0} \) can be retained; when \( \BF_{01}(y^{n}) \) is smaller than 1, the data provide evidence for \( \Hc_{1} \), the model that relaxes the general law. The larger the deviation from 1, the stronger the evidence. Importantly, in Jeffreys's framework the test-relevant parameter is fixed under \( \Hc_{0} \) and free to vary under \( \Hc_{1} \). The hypothesis \( \Hc_{0} \) is generally known as a `point-null' hypothesis. 
    
A perennial objection against point-null hypothesis testing---whether Bayesian or frequentist---is that in most practical applications, the point-null is never true exactly (e.g., \citealp{Bakan1966,Berkson1938,EdwardsEtAl1963,JonesTukey2000,KruschkeLiddell2018}; see also \citealp[p. 375]{Laplace17741986}). If this argument is accepted and \( \Hc_{0} \) is deemed to be false from the outset, then the test merely assesses whether or not the sample size was sufficiently large to detect the non-zero effect. This objection was forcefully made by Tukey:  

\begin{quotation}
``Statisticians classically asked the wrong question---and were willing to answer with a lie, one that was often a downright lie. They asked ``Are the effects of A and B different?'' and they were willing to answer ``no.'' 

All we know about the world teaches us that the effects of A and B are always different---in some decimal place---for any A and B. Thus asking ``Are the effects different?'' is foolish. \citep[p. 100]{Tukey1991}
\end{quotation}

This perennial objection has been rebutted in several ways (e.g., \citealp{Jeffreys1937Aristotelian,Jeffreys1961,KassRaftery1995}); in the current work we focus on the most common rebuttal, namely that the point-null hypothesis is a mathematically convenient approximation to a more realistic `peri-null' \citep{Tukey1995} hypothesis \( \Hc_{\widetilde{0}} \) that assigns the test-relevant parameter a distribution tightly concentrated around the value specified by the point-null hypothesis (e.g., \citealp[p. 416]{Good1967}; \citealp{BergerDelampady1987,Cornfield1966BayesianTest,Cornfield1969,Dickey1976,EdwardsEtAl1963,GeorgeMcCulloch1993,Jeffreys1935,Jeffreys1936Further,Gallistel2009,rousseau2007approximating,RouderEtAl2009Ttest}). For instance, in the case of the \( t \)-test the peri-null \( \Hc_{\widetilde{0}} \) could specify \( \delta \sim \pi(\delta \, | \, \Hc_{\widetilde{0}}) = \Nc (0, \kappa_{0}^{2}) \), where the width \( \kappa_{0} \) is set to a small value. 

Previous work has suggested that the approximation of a point-null hypothesis by an interval is reasonable when the width of that interval is half a standard error in width \citep{BergerDelampady1987,rousseau2007approximating} or one standard error in width \citep{Jeffreys1935}. Here we explore the consequences of replacing the point-null hypothesis \( \Hc_{0} \) by a peri-null hypothesis \( \Hc_{\widetilde{0}} \) from a different angle. We alter only the specification of the null-hypothesis \( \Hc_{0} \), which means that the alternative hypothesis \( \Hc_{1} \) now overlaps with \( \Hc_{\widetilde{0}} \). 

Below we show, first, that the effect on the Bayes factor of replacing \( \Hc_{0} \) with \( \Hc_{\widetilde{0}} \) is given by another Bayes factor, namely that between \( \Hc_{0} \) and \( \Hc_{\widetilde{0}} \) (cf. \citealp[p. 411]{MoreyRouder2011}). This `peri-null correction factor' is usually near 1, unless sample size grows large. For large sample sizes, we demonstrate that the Bayes factor for the peri-null \( \Hc_{\widetilde{0}} \) versus the alternative \( \Hc_{1} \) is bounded by the ratio of the prior ordinates evaluated at the maximum likelihood estimate. This proves earlier statements from \citet[pp. 411-412]{MoreyRouder2011} and confirms suggestions in \citet[p. 367]{Jeffreys1961} and \citet[p. 39, Eq. 2]{Jeffreys1973}. In other words, the Bayes factor for the peri-null hypothesis is inconsistent. 

Note that there exist several Bayes factor methods that have replaced point-null hypotheses with either peri-null hypotheses (e.g., Stochastic Search Variable Selection, \citealp{GeorgeMcCulloch1993}%
\footnote{\revise{``A similar setup in this context was considered by Mitchell and Beauchamp (1988), \nocite{MitchellBeauchamp1988} who instead used ``spike and slab'' mixtures. An important distinction of our approach is that we do not put a probability mass on $\beta_i = 0$.'' \citep[p. 883]{GeorgeMcCulloch1993}.}}) %
or with other hypotheses that have a continuous prior distribution close to zero (e.g., the sceptical prior proposed by \citealp{PawelHeld2020}). As far as evidence from the marginal likelihood is concerned, the results below show that these methods are inconsistent.

We end with a brief discussion on how a consistent method for hypothesis testing can be obtained without fully committing to a point-null hypothesis.  
	
\section{The Peri-Null Correction Factor}
Consider the three hypotheses discussed earlier: the point-null hypothesis \( \Hc_{0} \) fixes the test-relevant parameter to a fixed value (e.g., \( \delta = 0 \)); the peri-null hypothesis \( \Hc_{\widetilde{0}} \) assigns the test-relevant parameter a distribution that is tightly centered around the value of interest (e.g., \( \delta \sim \pi(\delta \, | \, \Hc_{\widetilde{0}}) = \Nc (0, \kappa_{0}^{2}) \) with \( \kappa_{0} \) small); and the alternative hypothesis \( \Hc_{1} \) assigns the test-relevant parameter a relatively wide prior distribution, \( \delta \sim \pi(\delta \, | \, \Hc_{1}) \). The Bayes factor of interest is between \( \Hc_{1} \) and \( \Hc_{\widetilde{0}} \), which can be expressed as the product of two Bayes factors involving \( \Hc_{0} \):
\begin{align}
\label{eq:BFcorrectionfactor}
\underbrace{ \frac{p(y^{n} \, | \, \Hc_{1})}{p(y^{n} \, | \, \Hc_{\widetilde{0}})}}_{\text{Peri-null } \BF_{1 \widetilde{0}} (y^{n}) } \quad = \underbrace{ \frac{p(y^{n} \, | \, \Hc_{1})}{p(y^{n} \, | \, \Hc_{0})}}_{\text{Point-null } \BF_{10} (y^{n}) } \quad \times \underbrace{ \frac{p( y^{n} \, | \, \Hc_{0})}{p(y^{n} \, | \, \Hc_{\widetilde{0}})}}_{\text{Correction factor } \BF_{0 \widetilde{0}} (y^{n}) } .
\end{align}
In words, the Bayes factor for the alternative hypothesis against the peri-null hypothesis equals the Bayes factor for the alternative hypothesis against the point-null hypothesis, multiplied by a correction factor (cf. \citealp{kass1992approximate}, \citealp{KassRaftery1995}, \citealp[p. 411]{MoreyRouder2011}). This correction factor quantifies the extent to which the point-null hypothesis outpredicts the peri-null hypothesis. With data sets of moderate size, and \( \kappa_{0} \) small, the peri-null and point-null hypotheses will make similar predictions, and consequently the correction factor will be close to 1. In such cases, the point-null can indeed be considered a mathematically convenient approximation to the peri-null.

\subsection{Example}
Consider the hypothesis that ``more desired objects are seen as closer'' \citep{BalcetisDunning2010}. In the authors' Study I, 90 participants had to estimate their distance to a bottle of water. Immediately prior to this task, 47 `thirsty' participants had consumed a serving of pretzels, whereas 43 `quenched' participants had drank as much as they wanted from four 8-oz glasses of water. In line with the authors' predictions, ``Thirsty participants perceived the water bottle as closer (\( M = 25.1 \) in., \( SD = 7.3 \)) than quenched participants did (\( M = 28.0 \) in., \( SD = 6.2 \))'' \citep[p. 148]{BalcetisDunning2010}, with \( t=2.00 \) and \( p=.049 \). A Bayesian point-null \( t \)-test concerning the test-relevant parameter \( \delta \) may contrast \( \Hc_{0}: \delta = 0 \) versus \( \Hc_{1} : \delta \in \R \) with a Cauchy distribution with location parameter 0 and scale \( \kappa_{1} \), the common default value \( \kappa_{1}=1/\sqrt{2} \) \citep{MoreyRouderBayesFactorPackage}. The resulting point-null Bayes factor is \( \BF_{10} = 1.259 \), a smidgen of evidence in favor of \( \Hc_{1} \). We may also compute a peri-null correction factor by contrasting \( \Hc_{0}: \delta = 0 \) against \( \Hc_{\widetilde{0}}: \delta \sim \Nc (0, \kappa_{0}^{2}) \), with \( \kappa_{0} = 0.01 \), say. The resulting peri-null correction factor%
\footnote{Calculated using the Summary Stats module in JASP, (e.g., \citealp{ly2018bayesian}, \url{jasp-stats.org}), and based on \citet{GronauEtAl2020Informed}.} %
is \( \BF_{0 \widetilde{0}} = 0.997 \), which means that, practically, it does not matter if the point-null or the peri-null is tested. With a larger value of \( \kappa_{0} = 0.05 \), we have \( \BF_{0\widetilde{0}} = 0.927 \), thus, a peri-null Bayes factor of \( \BF_{1 \widetilde{0}} = 1.167 \). The change from \( \BF_{1 0} = 1.259 \) to \( \BF_{1 \widetilde{0}} = 1.167 \) is utterly inconsequential.

The difference between the peri-null and point-null Bayes factor remains inconsequential for larger values of \( t \). When we change \( t=2.00 \) to \( t=4.00 \), the point-null Bayes factor equals \( \BF_{10} = 174 \), which according to Jeffreys's classification of evidence (e.g., \citealp[Appendix B]{Jeffreys1961}) is considered compelling evidence for \( \Hc_{1} \). With \( \kappa_{0} = 0.01 \), the peri-null correction factor equals \( \BF_{0\widetilde{0}} = 0.986 \) and consequently a peri-null Bayes factor equals of about \( 172 \) in favor of \( \Hc_{1} \) over \( \Hc_{\widetilde{0}} \). With \( \kappa_{0} = 0.05 \), the peri-null correction factor equals \( \BF_{0 \widetilde{0}} = 0.713 \) and \( \BF_{1\widetilde{0}} \approx 124 \). In absolute numbers, the change from 174 to 124 may appear considerable, but with equal prior model probabilities this translates to a modest difference in posterior probabilities: \( P(\Hc_{1} \, | \, y^{n}) = 174/175 \approx 0.994 \) versus \( 124/125 = 0.992 \). %

The peri-null correction factor does become influential when sample size is large. As we prove in the next section, the peri-null Bayes factor is inconsistent and converges to the ratio of prior ordinates under \( \Hc_{1} \) and \( \Hc_{\widetilde{0}} \) at the maximum likelihood estimate. 

\section{The Peri-Null Bayes Factor is Inconsistent}
Historically, the main motivation for the development of the Bayes factor was the desire to be able to obtain arbitrarily large evidence for a general law: ``We are looking for a system that will in suitable cases attach probabilities near 1 to a law.'' (\citealp[p. 88]{Jeffreys1977}; see also \citealp{EtzWagenmakers2017,ly2020bayesian,WrinchJeffreys1921}). 

Statistically, this desideratum means that we want Bayes factors to be consistent, which implies that, as sample size increases, (i) \( \BF_{10}(Y^{n}) \) tends to zero when the data are generated under the null model, whereas (ii) \( \BF_{01}(Y^{n}) \) tends to zero when the data are generated under the alternative model \( \Hc_{1} \), that is, %
\begin{align} 
 \BF_{10} (Y^{n}) \overset{\Pf_{\theta}}{\rightarrow} 0 \textnormal{ if } \Pf_{\theta} \in \Hc_{0} , \textnormal{ and }  \BF_{01} (Y^{n}) \overset{\Pf_{\theta}}{\rightarrow} 0 \textnormal{ if } \Pf_{\theta} \in \Hc_{1} .
\end{align} 
Thus, regardless of the chosen prior model probabilities \( P ( \Hc_{0} ), P ( \Hc_{1} ) \in (0, 1) \), 
\begin{align} 
\label{eqBfConsistency}
P( \Hc_{0} \, | \, Y^{n}) \overset{\Pf_{\theta}}{\rightarrow} 1 \textnormal{ if } \Pf_{\theta} \in \Hc_{0}, \textnormal{ and }  P( \Hc_{1} \, | \, Y^{n}) \overset{\Pf_{\theta}}{\rightarrow} 1 \textnormal{ if } \Pf_{\theta} \in \Hc_{1}, 
\end{align} 
where \( \Pf_{\theta} \) refers to the data generating distribution, here, \( Y_{i} \iidSim \Pf_{\theta} \), and where  \( X_{n} \overset{\Pf_{\theta}}{\rightarrow} X \) denotes convergence in probability, that is, \( \lim_{n \rightarrow \infty} \Pf_{\theta} ( | X_{n} - X | > \epsilon ) = 0  \) as usual. 

Below we prove that the peri-null Bayes factor \( \BF_{1 \widetilde{0}}(Y^{n}) \) is inconsistent (cf. suggestions by \citealp[p. 367]{Jeffreys1961}; \citealp[p. 39, Eq. 2]{Jeffreys1973}; and the statements by \citealp[p. 411-412]{MoreyRouder2011}). The proof relies on the observation that the replacement of the point-null restriction on the test-relevant parameter, i.e., \( \Hc_{0} : \delta = 0 \), where \( \theta = (\delta, \theta_{0}) \) as before, yields a peri-null model that defines the same likelihood function as the alternative model. Consequently, the numerator and the denominator of the peri-null Bayes factor \( \BF_{1 \widetilde{0}}(Y^{n}) \) only differ in how the priors are specified. 

The inconsistency of peri-null Bayes factors then follows quite directly from Laplace's method \citep{Laplace17741986} and consistency of the maximum likelihood estimator (MLE). Both Laplace's method and consistency of the MLE hold under weaker conditions than stated here, namely, for absolute continuous priors (e.g., \citealp[Chapter~{10}]{van1998asymptotic}), and regular parametric models (e.g., \citealp[Chapter~{7}]{van1998asymptotic}; \citealp[Appendix~{E}]{ly2017tutorial}). These models only need to be one time differentiable with respect to \( \theta \) in quadratic mean and have non-degenerate Fisher information matrices that are continuous in \( \theta \) with determinants that are bounded away from zero and infinity. The inconsistency of the peri-null Bayes factor is therefore expected to hold more generally. 

We show that under the stronger conditions of \citet{kass1990validity}, the asymptotic sampling distribution of peri-null Bayes factors can be easily derived. These stronger conditions imply that the model is regular for which we know that the MLE is not only consistent, but also locally asymptotically normal with a variance equal to the inverse observed Fisher information matrix at \( \hat{\theta} \) with entries
\begin{align} 
[ \hat{I}(\hat{\theta})]^{a, b} = - \frac{1}{n} \sum_{i=1}^{n} \Big ( \tfrac{ \partial^{2}}{ \partial \theta_{a} \partial \theta_{b}} \log f(Y_{i} \, | \, \theta) \Big ) \Bigg |_{\theta = \hat{\theta}},
\end{align} 
see for instance \citet{ly2017tutorial} for details. 

\begin{theorem}[Limit of a peri-null Bayes factor] 
\label{thmLimitPerinullBf}
Let \( Y^{n} = ( Y_{1}, \ldots, Y_{n}) \) with \( Y_{i} \iidSim \Pf_{\theta} \in \Pc_{\Theta} \), where \( \Pc_{\Theta} \) is an identifiable family of distributions that is Laplace-regular \citep{kass1990validity}. This implies that \( \Pc_{\Theta} \) admits densities \( f(y^{n} \, | \, \theta) \) with respect to the Lebesgue measure that are six times continuously differentiable in \( \theta \) at the data-governing parameter \( \theta \in \Theta \subset \R^{p} \) and \( \Theta \) open with non-empty interior. Furthermore, assume that the (peri-null) prior densities \( \pi( \theta \, | \, \Hc_{\widetilde{0}}) \) and \( \pi( \theta \, | \, \Hc_{1}) \) assign positive mass to a neighborhood at the data-governing parameter \( \theta \) and are four times continuously differentiable at \( \theta \); then \( \BF_{1 \widetilde{0}}(Y^{n}) \overset{\Pf_{\theta}}{\rightarrow} \tfrac{ \pi( \theta \, | \, \Hc_{1})}{\pi( \theta \, | \, \Hc_{\widetilde{0}})} \). \( \hfill \diamond \)
\end{theorem} 

\begin{proof} 
The condition that the model is Laplace-regular allows us to employ the Laplace method to approximate the numerator and the denominator of the peri-null Bayes factor by
\begin{align}
\label{eqLaplaceApprox}
p(Y^{n} \, | \, \Hc_{j}) & = f(Y^{n} \, | \, \hat{\theta}) \big ( \tfrac{2 \pi}{n} \big )^{\frac{p}{2}} | \hat{I}(\hat{\theta}) |^{-\frac{1}{2}}  \pi(\hat{\theta} \, | \, \Hc_{j}) \\
\nonumber
& \times \Big ( 1 + \tfrac{ C^{1}(\hat{\theta} \, | \, \Hc_{j}) }{n} + \tfrac{ C^{2}(\hat{\theta} \, | \, \Hc_{j}) }{n^{2}} + \Oc(n^{-3}) \Big ) ,
\end{align}
where \( C^{1}(\hat{\theta} \, | \, \Hc_{j}) \) and \( C^{2}(\hat{\theta} \, | \, \Hc_{j}) \) for \( j= \widetilde{0}, 1 \) are bounded error terms of the Laplace approximation (cf. \citealp{kass1990validity}) and given explicitly by \refEq{laplaceOrder1} and \refEq{laplaceOrder2} based on the notation introduced in \refApp{appLaplace}.

From the fact that the peri-null and the alternative models define the same likelihood function, thus, have the same maximum likelihood estimator, and only differ in how the priors concentrate on the parameters, we conclude that 
\begin{align} 
\label{eqPeriNullBf} 
\BF_{1 \widetilde{0}}(Y^{n}) = \frac{ \pi(\hat{\theta} \, | \, \Hc_{1}) \Big [1 + \tfrac{C^{1}(\hat{\theta} \, | \, \Hc_{1})}{n}  + \Oc(n^{-2}) \Big ] }{\pi(\hat{\theta} \, | \, \Hc_{\widetilde{0}}) \Big [1 + \tfrac{C^{1}(\hat{\theta} \, | \, \Hc_{\widetilde{0}})}{n} + \Oc(n^{-2}) \Big ] } . 
\end{align}
Identifiability and the regularity conditions on the model imply that the maximum likelihood estimator is consistent, thus, \( \hat{\theta} \overset{\Pf_{\theta}}{\rightarrow} \theta \) (e.g., \citealp[Chapter~{5}]{van1998asymptotic}). As all functions of \( \hat{\theta} \) in \refEq{eqPeriNullBf} are smooth at \( \theta \), the continuous mapping theorem applies and the assertion follows. \( \hfill \qed \)
\end{proof} 

\refThm{thmLimitPerinullBf} implies that \( \BF_{1 \widetilde{0}}(Y^{n}) \) is inconsistent; for all data-governing parameter values with a neighborhood that receives positive mass from both priors, the peri-null Bayes factor approaches a limit that is given by the ratio of prior densities evaluated at the data-governing \( \theta \) as \( n \) increases. Note that this holds in particular for the test point of interest, e.g., \( \delta = 0 \), which has a neighborhood that the peri-null prior assigns positive mass to. This inconsistency result can be intuited as follows. The peri-null Bayes factor compares two marginal likelihoods with the same data-distribution (or sampling) model, but different prior distributions on the same parameter space; hence, the peri-null Bayes factor effectively assesses which prior performs best, and this should not matter asymptotically (i.e., as the data accumulate, the posterior distributions of the two models converge, and consequently the change in the Bayes factor will converge as well).%
\footnote{We thank the first anonymous reviewer for providing this intuition.} %

From \refThm{thmLimitPerinullBf} it follows that the Bayes factor comparing the alternative \( \Hc_{1} : \delta \in \Delta = \R \) against a directed hypothesis, say, \( \Hc_{+} : \delta > 0 \), is also inconsistent. For data under any \( \delta > 0 \), the associated Bayes factor \( \BF_{1+}(Y^{n}) \) then converges in probability to \( \Pi_{u}( \{ \delta > 0 \} )/ \Pi_{u}(\Delta) \), where \( \Pi_{u}(B) = \int_{B} \pi_{u}(\theta) \der \theta \) with \( \pi_{u} \) the unnormalized prior on \( \delta \). 

The limit described in \refThm{thmLimitPerinullBf} can also be derived differently. For instance, Theorem 1 (ii) of \citet{dawid2011posterior} can be applied twice: once to approximate the logarithm of the marginal likelihood of the alternative model, and once for the null model.%
\footnote{We thank the second anonymous reviewer for bringing this reference to our attention. When comparing our \refThm{thmLimitPerinullBf} to that of Theorem 1 (ii) of \citet{dawid2011posterior}, it is worth noting that the apparent difference in the order of the remainder term vanishes once the MLE \( \hat{\theta} \) in \refEq{eqLaplaceApprox} is replaced by \( \hat{\theta} = \theta + h / \sqrt{n} \), which holds in probability for large \( n \) for regular models.} %
Another way to derive the limit in \refThm{thmLimitPerinullBf} is by using the generalized Savage-Dickey density ratio \citep{VerdinelliWasserman1995} and by exploiting the transitivity of the Bayes factor. \refThm{thmLimitPerinullBf}, however, can be more straightforwardly extended to characterize the asymptotic behavior of the peri-null Bayes factor. 

The limiting value of the peri-null Bayes factor is not representative when \( n \) is small or moderate. \refThm{thmAsymptoticDistributionPeriNull} below shows that the sampling mean of \( \log \BF_{1 \widetilde{0}}(Y^{n}) \) is expected to be of smaller magnitude than its limiting value. In other words, the limit in \refThm{thmLimitPerinullBf} should be viewed as an upper bound under the alternative and a lower bound under the null. 

This theorem exploits the fact that without a point-null hypothesis the gradients of the densities \( \pi(\theta \, | \, \Hc_{1}) \) and \( \pi(\theta \, | \, \Hc_{\widetilde{0}}) \) are of the same dimension, which implies that the gradient  \( \tfrac{ \partial}{\partial \theta } \log \big ( \tfrac{\pi ( \theta \, | \, \Hc_{1})}{\pi( \theta \, | \, \Hc_{\widetilde{0}})} \big ) \) is well-defined. As such, the delta method can be used to show that the peri-null Bayes factor inherits the asymptotic normality property of the MLE. 

To state the theorem we write \( D \) for the differential operator with respect to \( \theta \), e.g., \( [D^{1} \pi(\theta \, | \, \Hc_{j})] = \frac{ \partial}{\partial \theta} \pi(\theta \, | \, \Hc_{j}) \) denotes the gradient, and \( [ D^{2} \pi(\theta \, | \, \Hc_{j})]  = [ \frac{ \partial^{2}}{\partial \theta \partial \theta} \pi(\theta \, | \, \Hc_{j}) ]\) denotes the Hessian matrix. 

\begin{theorem}[Asymptotic sampling distribution of a peri-null Bayes factor]
\label{thmAsymptoticDistributionPeriNull}
Under the regularity conditions stated in \refThm{thmLimitPerinullBf} and for all data-governing parameters \( \theta \) for which %
\begin{align}
\label{eqNoZeroLogPriorRatioGradient}
\dot{v}(\theta) :=  [ D^{1} \log \big ( \tfrac{\pi( \theta \, | \, \Hc_{1})}{\pi( \theta \, | \, \Hc_{\widetilde{0}})} \big )] \neq 0 \in \R^{p}, 
\end{align}
the asymptotic sampling distribution of the logarithm of the peri-null Bayes factor is normal, that is, %
\begin{align}
\sqrt{n} \Big ( \log \BF_{1 \widetilde{0}}(Y^{n}) - \log \big ( \tfrac{\pi ( \theta \, | \, \Hc_{1})}{\pi( \theta \, | \, \Hc_{\widetilde{0}})} \big ) - E(\theta, n) \Big ) \overset{\Pf_{\theta}}{\rightsquigarrow} \Nc \left ( 0 , \dot{v}(\theta)^{T} I^{-1}(\theta) \dot{v}(\theta) \right ) ,
\end{align}
where \( \overset{\Pf_{\theta}}{\rightsquigarrow} \) denotes convergence in distribution under \( \Pf_{\theta} \) and where %
\begin{align}
E(\theta, n) &= \log \big ( \tfrac{ 1 + C^{1}(\theta \, | \, \Hc_{1}) /n + C^{2}(\theta \, | \, \Hc_{1}) /n^{2} }{ 1 + C^{1}(\theta \, | \, \Hc_{\widetilde{0}}) /n + C^{2}(\theta \, | \, \Hc_{\widetilde{0}}) /n^{2} } \big ) , 
\end{align}
is a bias term that is asymptotically negligible, and where \( C^{1}(\theta \, | \, \Hc_{j}) \) and \( C^{2}(\theta \, | \, \Hc_{j}) \) are given explicitly by \refEq{laplaceOrder1} and \refEq{laplaceOrder2} based on the notation introduced in \refApp{appLaplace}.

For all \( \theta \) for which \( \dot{v}(\theta) = 0 \), but \( \ddot{v}(\theta) :=  [ D^{2} \log \big ( \tfrac{\pi( \theta \, | \, \Hc_{1})}{\pi( \theta \, | \, \Hc_{\widetilde{0}})} \big )] \neq 0  \in \R^{p \times p} \), the asymptotic distribution of \( \log \BF_{1 \widetilde{0}}(Y^{n}) \) has a quadratic form, that is,
\begin{align}
n \Big ( \log \BF_{1 \widetilde{0}}(Y^{n}) - \log \big ( \tfrac{\pi ( \theta \, | \, \Hc_{1})}{\pi( \theta \, | \, \Hc_{\widetilde{0}})} \big ) - E(\theta, n) \Big ) \overset{\Pf_{\theta}}{\rightsquigarrow} Z^{T} I^{-1/2}(\theta) \ddot{v}(\theta) I^{-1/2}(\theta) Z, 
\end{align}
where \( Z \sim \Nc( 0 , I ) \) with \( I \in \R^{p \times p} \) the identity matrix. \( \hfill \diamond \)
\end{theorem} 

\begin{proof} 
The proof depends on (another) Taylor series expansion, see \refApp{appLaplace} for full details. Firstly, we recall that \( \sqrt{n} ( \hat{\theta} - \theta) \overset{\theta}{\rightsquigarrow} \Nc ( 0, I^{-1}(\theta)) \). To relate this asymptotic distribution to that of \( \log \BF_{1 \widetilde{0}}(Y^{n}) \), we note that \refEq{eqPeriNullBf} is, up to a decreasing error in \( n \), a smooth function of the maximum likelihood estimator. The goal is to ensure that the error terms \( 1 + C^{1}( \theta \, | \, \Hc_{j})/n + C^{2}( \theta \, | \, \Hc_{j})/n^{2} \) are asymptotically negligible. A Taylor series expansion at the data-governing \( \theta \) shows that %
\begin{align}
\label{eqBfOrder0}
\log \BF_{1 \widetilde{0}}(Y^{n}) & = \log \big ( \tfrac{\pi ( \theta \, | \, \Hc_{1})}{\pi( \theta \, | \, \Hc_{\widetilde{0}})}  \big ) + E (\theta, n) \\
\nonumber
& + (\hat{\theta} - \theta)^{T} \big ( \dot{v}(\theta) +  [D^{1} E (\theta, n)]   \big ) \\
\nonumber
& +  ( \hat{\theta} - \theta)^{T} \tfrac{  ( \ddot{v}(\theta) + [ D^{2} E (\theta, n)  ] )}{2} (\hat{\theta} - \theta)+ \Oc_{P}(n^{-3/2}) .
\end{align}
The asymptotic normality result follows after rearranging \refEq{eqBfOrder0}, a multiplication of \( \sqrt{n} \) on both sides, and an application of Slutsky's lemma. 

To conclude that the bias term \( E(\theta, n) \) is indeed asymptotically negligible, note that \( \log ( 1 + x/n) \approx x/n \) as \( n \rightarrow \infty \) and therefore \( D^{k} E(\theta, n) = \Oc \Big ( \frac{1}{n} D^{k} \big \{ C^{1}(\theta \, | \, \Hc_{1}) - C^{1}(\theta \, | \, \Hc_{\widetilde{0}}) \big \}  \Big ) \) for all \( k \leq 3 \). The approximation \( \log ( 1 + x/n) \approx x/n \) requires \( C^{k}(\theta \, | \, \Hc_{j})  \) for \( k=1,2 \) and \( j=\widetilde{0}, 1 \) to be of similar magnitude, but this is typically not the case when \( \kappa_{0} \) is relatively small compared to \( \kappa_{1} \). The bias is, therefore, expected to decay much more slowly. 

Similarly, when \( \dot{v}(\theta) \) is zero, but \( \ddot{v}(\theta) \) not, we have
\begin{align}
n \log \BF_{1 \widetilde{0}}(Y^{n}) & = n \big (  \log \big ( \tfrac{\pi ( \theta \, | \, \Hc_{1})}{\pi( \theta \, | \, \Hc_{\widetilde{0}})}  \big ) +  E (\theta, n)) \big ) \\
\nonumber
& +   \sqrt{n} ( \hat{\theta} - \theta)^{T} \tfrac{  ( \ddot{v}(\theta) + \Oc(n^{-1})  ] )}{2} \sqrt{n} (\hat{\theta} - \theta)+ \Oc_{P}(n^{-1/2}) .
\end{align}
Since \( \sqrt{n} ( \hat{\theta} - \theta) \overset{\Pf_{\theta}}{\rightsquigarrow} \Nc ( 0, I(\theta)^{-1}) \), the second order result follows. \( \hfill \qed \)
\end{proof}

\refThm{thmAsymptoticDistributionPeriNull} also shows that under the alternative hypothesis, \( \log \BF_{1 \widetilde{0}}(Y^{n}) \) is expected to increase towards the limiting value \( \log \big ( \tfrac{\pi ( \theta \, | \, \Hc_{1})}{\pi( \theta \, | \, \Hc_{\widetilde{0}})} \big ) \) as \( n \rightarrow \infty  \) whenever \( E(\theta, n) < 0 \). The bias is expected to be negative, because if the data-governing parameter \( \delta \) is far from zero, but the peri-null prior is specified such that it is peaked at zero, the Laplace approximations become less accurate. In other words, for fixed \( n \) and \( \delta \neq 0  \), we typically have \( C^{1}(\theta \, | \, \Hc_{1})  \leq C^{1}(\theta \, | \, \Hc_{\widetilde{0}})  \) and \( C^{2}(\theta \, | \, \Hc_{1}) \leq C^{2}(\theta \, | \, \Hc_{\widetilde{0}})  \) and, therefore, \( E(\theta, n) < 0 \). This intuition can be made rigorous using the explicit formulas for \( E(\theta, n) \) provided by \refEq{laplaceOrder1} and \refEq{laplaceOrder2} from \refApp{appLaplace}, as is shown in the following example. 

\section{Example}
\begin{figure}[h]
\begin{tabular}{cc}
    \begin{minipage}{.5 \textwidth}
    \centering
    \qquad \scalebox{1.25}{\( \kappa_{0}=0.05 \)}
    \includegraphics[width= \linewidth]{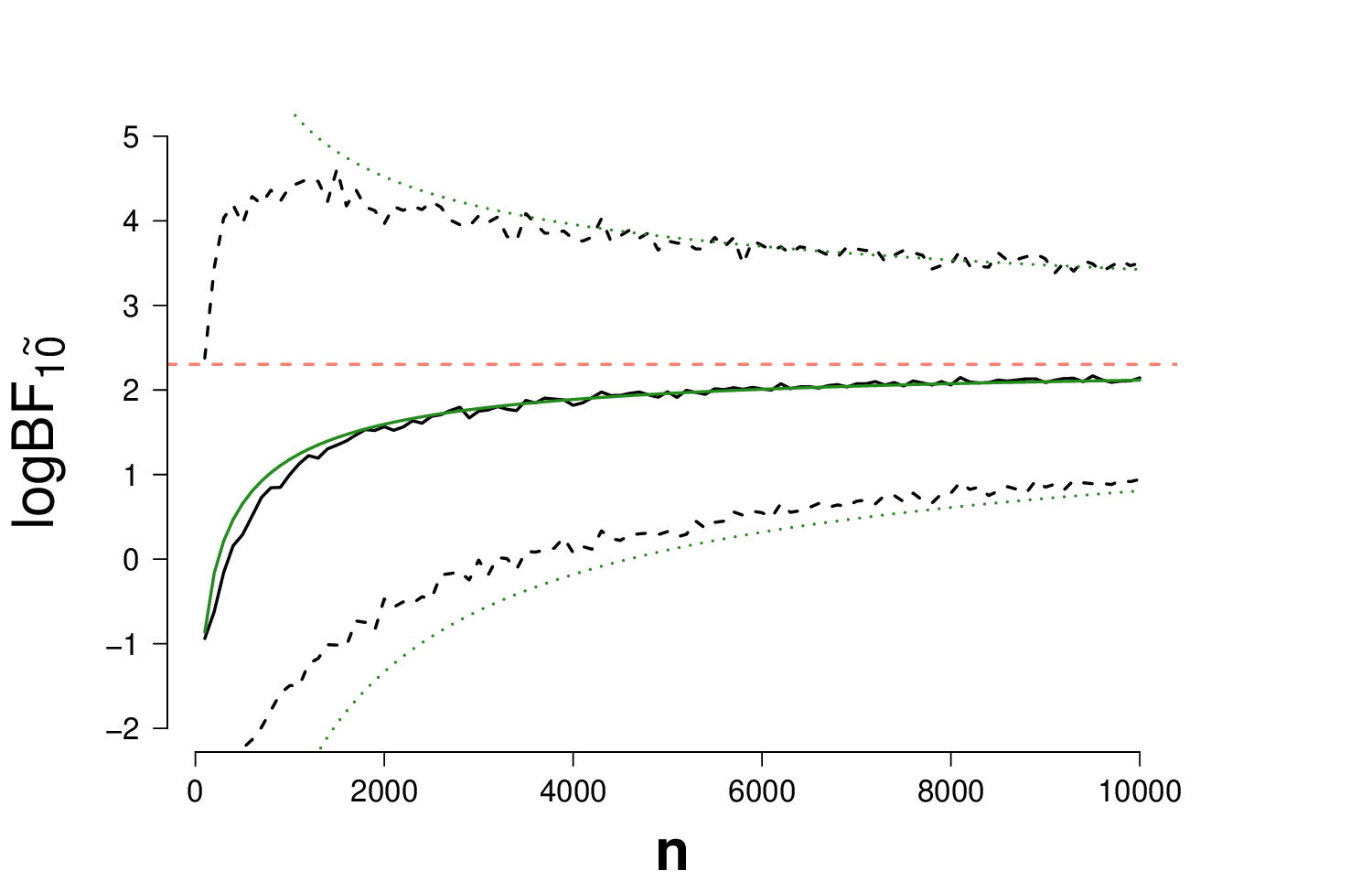}
    \end{minipage}  & %
    \begin{minipage}{.5 \textwidth}
    \centering
    \qquad \scalebox{1.25}{\( \kappa_{0}=0.10 \)}
    \includegraphics[width=\linewidth]{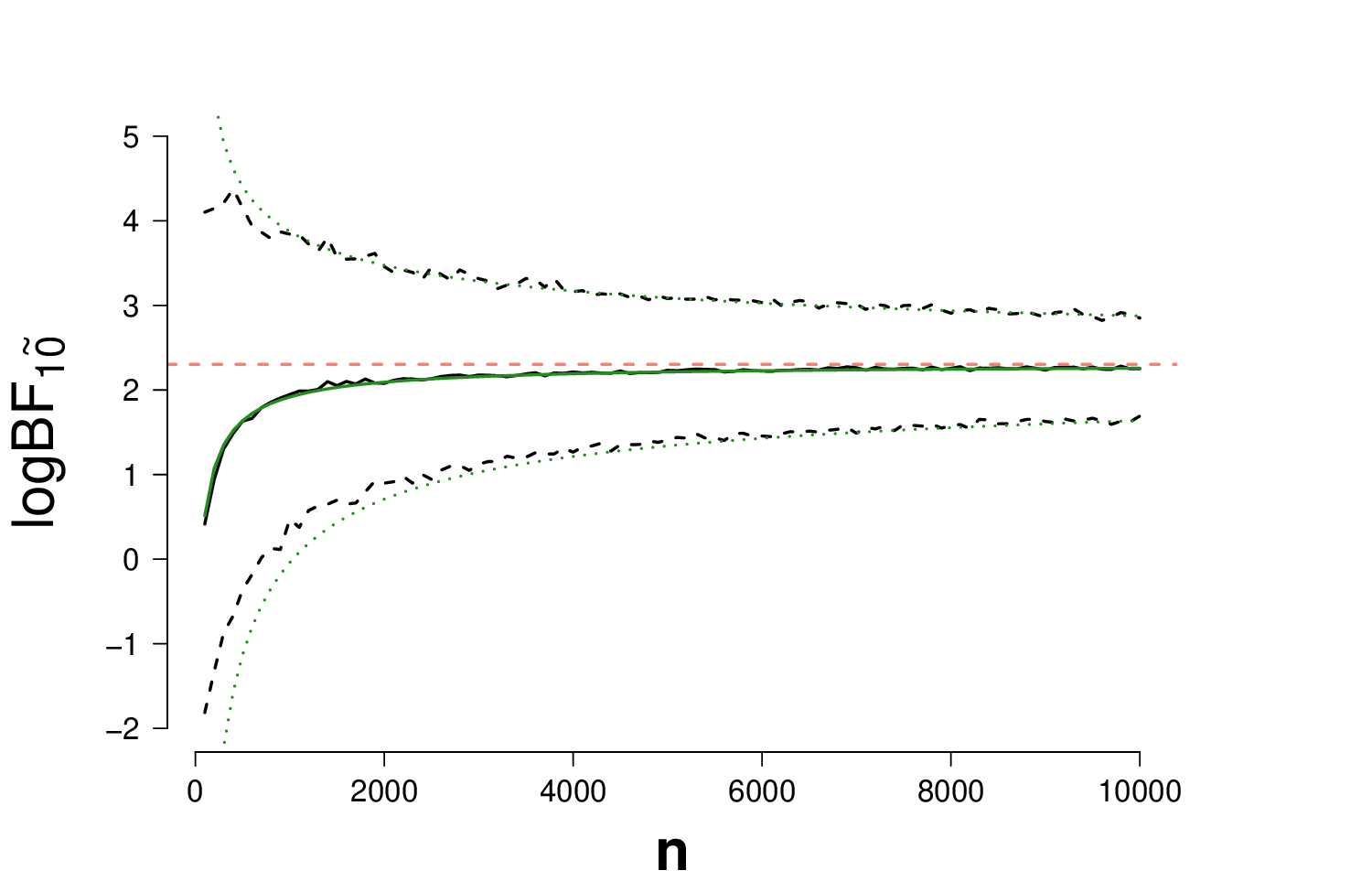}
    \end{minipage}  %
    \\
    \begin{minipage}{.5 \textwidth}
    \centering
    \includegraphics[width= \linewidth]{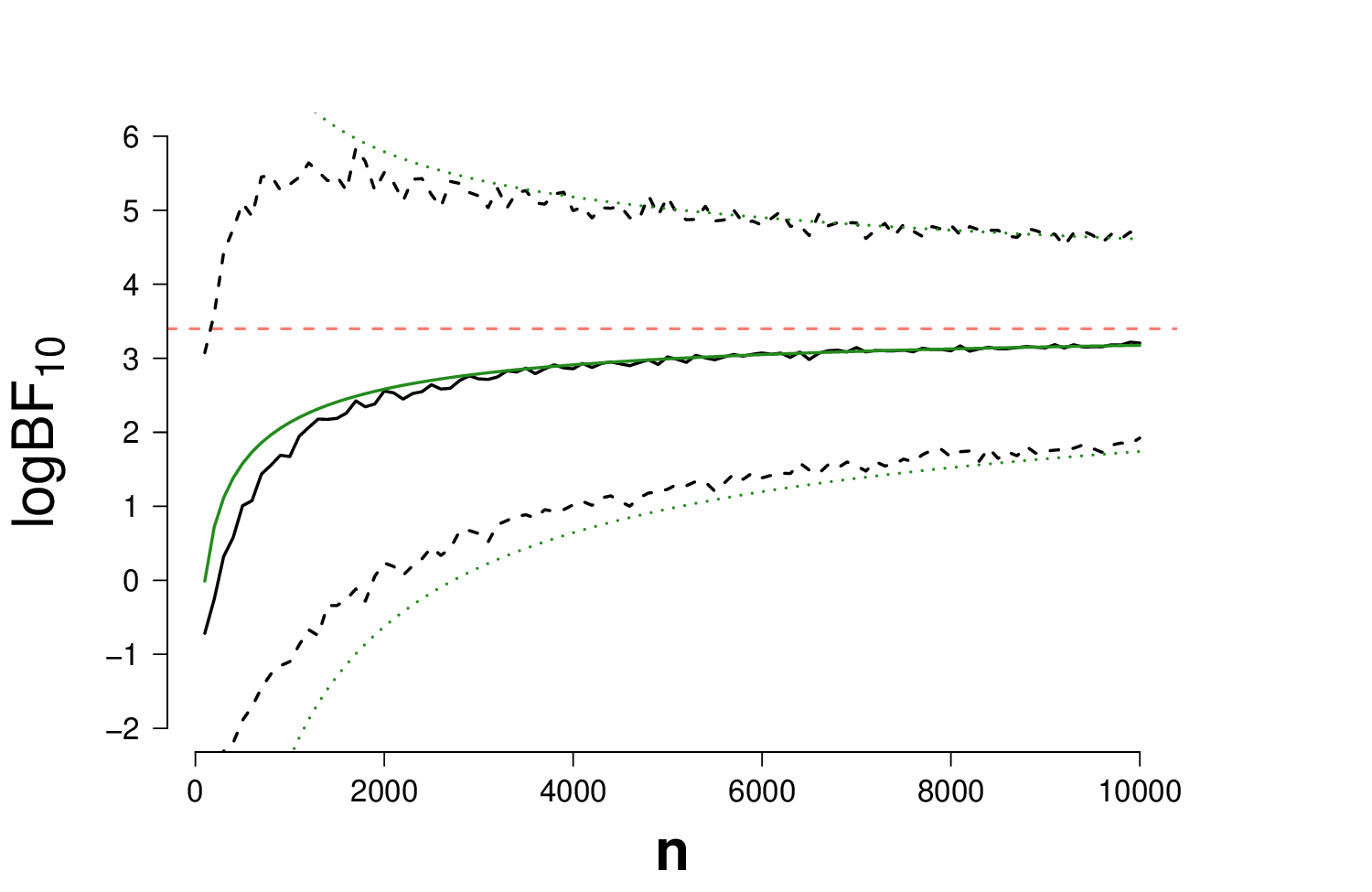}
    \end{minipage}  & %
    \begin{minipage}{.5 \textwidth}
    \centering
    \includegraphics[width=\linewidth]{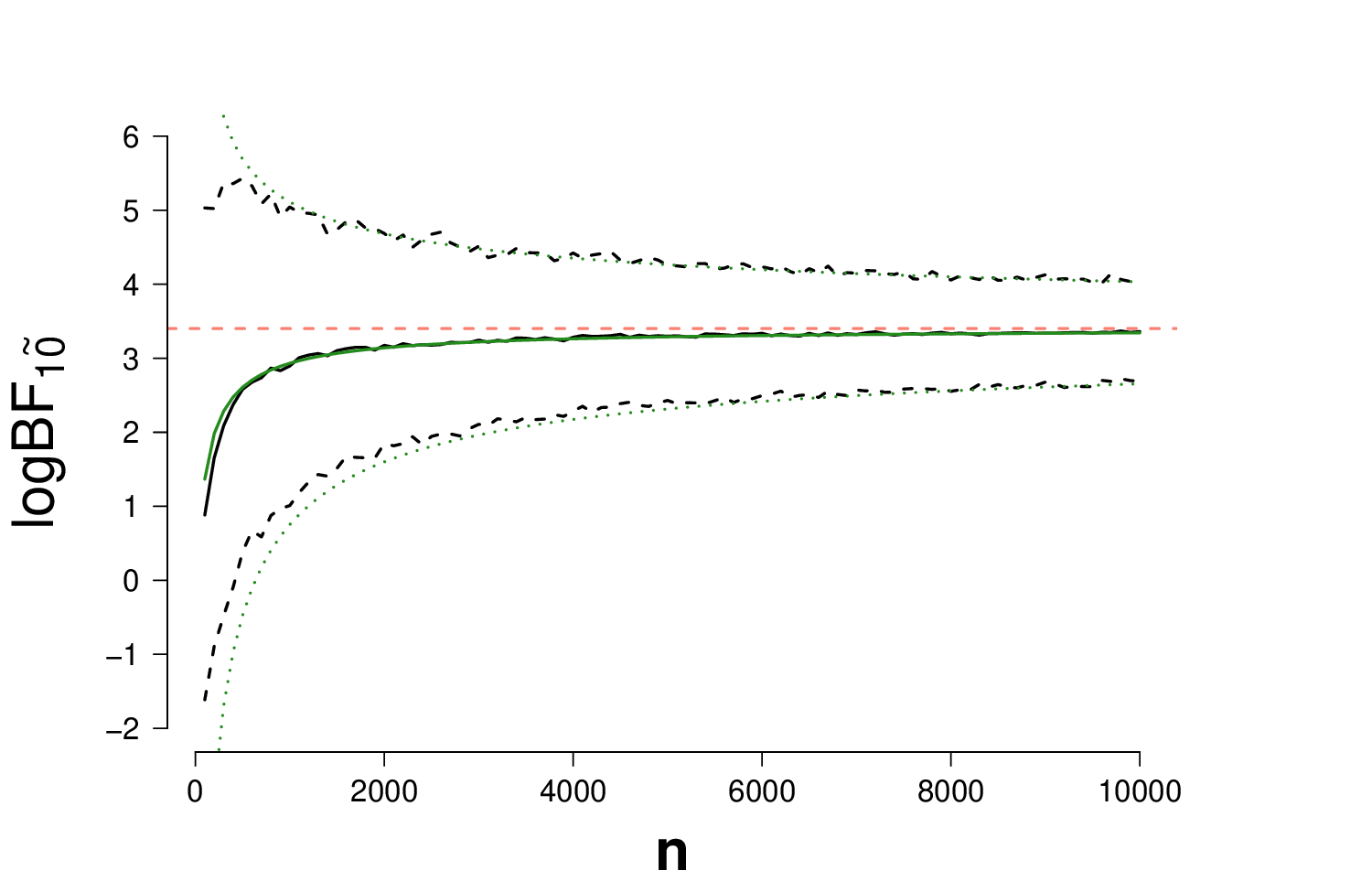}
    \end{minipage}  %
\end{tabular}
\caption{Under the alternative, the logarithm of the peri-null Bayes factor \( t \)-test is asymptotically normal with a mean (i.e., the solid curves) that increases to the limit, e.g., \( \log \BF_{1 \widetilde{0}} = \log ( 10) \) and \( \log \BF_{1 \widetilde{0}} = \log ( 30) \) in the top and bottom row respectively. The black and green curves correspond to the simulated and asymptotic normal sampling distribution respectively. The dotted curves show the 97.5\% and 2.5\% quantiles of the respective sampling distribution. Note that the convergence to the upper bound is slower when the peri-null is more concentrated, e.g., compare the left to the right column.}
\label{figBfCNAlternative}
\end{figure}

We consider a Bayesian \( t \)-test and for the peri-null Bayes factor use the priors 
\begin{align}
\label{eqPeriNullExample}
\pi(\delta, \sigma \, | \, \Hc_{1}) \propto \textnormal{Cauchy}(\delta \, ; \, 0, \kappa_{1}) \sigma^{-1} \text{ and } \pi(\delta, \sigma \, | \, \Hc_{\widetilde{0}}) \propto \Nc ( \delta \, ; \, 0, \kappa_{0}^2) \sigma^{-1}.
\end{align}
Note that \( \pi(\delta, \sigma \, | \, \Hc_{1}) \) is chosen as in the default Bayesian \( t \)-test \citep{jeffreys1948theory2, ly2016harold, ly2016evaluation, RouderEtAl2009Ttest}, where \( \kappa_{1} > 0 \) denotes the scale parameter of the Cauchy distribution on standardized effect size \( \delta = \mu / \sigma \), and \( \sigma \propto \sigma^{-1} \) implies that the standard deviation common in both models is proportional to \( \sigma^{-1} \) (for advantages of this choice see \citealp{hendriksen2021optional, grunwald2019safe}). For data-governing parameters \( \theta = (\mu, \sigma) \), where \( \mu \) is the population mean, \refThm{thmLimitPerinullBf} shows that as \( n \rightarrow \infty \)
\begin{align}
\label{eqPeriNullLimit}
\log  \BF_{1 \widetilde{0}}(Y^{n} \, ; \, \kappa_{0}, \kappa_{1}) & \overset{\Pf_{\theta}}{\rightarrow} \log \left ( \frac{ \sqrt{2} \kappa_{0} \exp (  \frac{ \mu^{2}}{2 \kappa_{0}^{2} \sigma^{2}})}{\sqrt{\pi} \kappa_{1} ( 1 + \big [\frac{ \mu}{\kappa_{1} \sigma} \big ]^{2})} \right ) =: v(\theta) .
\end{align}
Direct calculations show that \( \dot{v}(\theta) = 0 \) only when \( \mu=0 \). Hence, under the alternative \( \mu \neq 0 \), the logarithm of these peri-null Bayes factor \( t \)-tests are asymptotically normal with an approximate variance of
\begin{align}
\label{eqExAVar}
\frac{ (\mu^{4} + 2 \mu^{2} \sigma^{2}) (\mu^{2} + (\kappa_{1}^{2} - 2 \kappa_{0}^{2}) \sigma^{2})^{2} }{2 \kappa_{0}^{4} \sigma^{4} ( \mu^{2} + \kappa_{1}^{2} \sigma^{2})^{2} n} .
\end{align}
To characterize the asymptotic mean we also require the bias term \( E(\theta, n) \). An application of \refEq{laplaceOrder1} and \refEq{laplaceOrder2} from \refApp{appLaplace} shows that for the problem at hand the bias term comprises of 
\begin{align} %
C^{1}(\theta \, | \, \Hc_{1}) & = \tfrac{ 13 \mu^{4} + ( 18 + 2 \kappa_{1}^{2}) \sigma^{2} \mu^{2} + ( \kappa_{1}^{2} - 6) \kappa_{1}^{2} \sigma^{4}}{6 ( \mu^{2} + \kappa_{1}^{2} \sigma^{2})^{2}} , \\ %
C^{2}(\theta \, | \, \Hc_{1}) & = \tfrac{ 780 \mu^{6} + ( 1110 + 3127 \kappa_{1}^{2}) \sigma^{2} \mu^{4} + ( 6020 + 4462 \kappa_{1}^{2} ) \kappa_{1}^{2} \sigma^{4} \mu^{2} +  ( 5091 \kappa_{1}^{2} - 1426) \kappa_{1}^{4} \sigma^{6}}{-96(\mu^{2} +\kappa_{1}^{2} \sigma^{2})^{3}} , \\ %
C^{1}(\theta \, | \, \Hc_{\widetilde{0}}) & = \tfrac{3 \mu^{4} + 6 \sigma^{2} \mu^{2} + \kappa_{0}^{2} \sigma^{4} ( 2 \kappa_{0}^{2} - 6)}{12 \kappa_{0}^{4} \sigma^{4}},  \\ %
C^{2}(\theta \, | \, \Hc_{\widetilde{0}}) & = \tfrac{ 124 \mu^{6} + ( 264 - 2369 \kappa_{0}^{2}) \sigma^{2} \mu^{4} + (10811 \kappa_{0}^{2} - 2218) \kappa_{0}^{2} \sigma^{4} \mu^{2} + 2(713 - 5091 \kappa_{0}^{2}) \kappa_{0}^{4} \sigma^{6}}{192 \kappa_{0}^{6} \sigma^{6}} . %
\end{align}
More concretely, under \( \mu = 0.167 \) and \( \sigma = 1 \), \( \log  \BF_{1 \widetilde{0}}(Y^{n} \, ; \, 0.05, 1) \) converges in probability to \( \log(10) \). This limit is depicted as the pink dashed horizontal line in the top left subplot of \refFig{figBfCNAlternative}. 

\begin{sloppypar}
This subplot also shows the mean (solid green curve) and the 97.5\% and 2.5\% quantiles (dotted green curves above and below the solid curve respectively) based on the asymptotic normal result of \refThm{thmAsymptoticDistributionPeriNull}. The black curves represent the analogous quantities based on simulated normal data with \( \mu=0.167 \), \( \sigma=1 \) based on 1,000 replications at sample sizes \( n=100, 200, 300, \ldots, 10000 \). 
\end{sloppypar}

Observe that for small sample sizes, the simulated peri-null Bayes factors are more concentrated on small values. In this regime the concentration of the peri-null prior dominates, and the Laplace approximation of \( p(Y^{n} \, | \, \Hc_{\widetilde{0}}) \) is still inaccurate. 

As expected, the Laplace approximation becomes accurate sooner, whenever the peri-null prior is less concentrated. The top right subplot depicts results of \( \log  \BF_{1 \widetilde{0}}(Y^{n} \, ; \, 0.10, 1)  \) under \( \mu=0.314 \) and \( \sigma=1 \), which converges in probability to \( \log(10) \). 

Similarly, the asymptotic normal distribution becomes adequate at a smaller sample size for larger population means \( \mu \). The bottom left subplot corresponds to \( \log  \BF_{1 \widetilde{0}}(Y^{n} \, ; \, 0.05, 1)  \) under \( \mu=0.182 \) and \( \sigma=1 \), whereas the bottom right subplot corresponds to \( \log  \BF_{1 \widetilde{0}}(Y^{n} \, ; \, 0.10, 1)  \) under \( \mu=0.348 \) and \( \sigma=1 \). The logarithms of both peri-null Bayes factors converge in probability to \( \log(30) \). 

In sum, the plots show that under the alternative hypothesis the asymptotic normal distribution approximates the sampling distribution of the logarithm of the peri-null Bayes factor quite well, and it approximates better when the peri-null prior is less concentrated.

Under the null hypothesis \( \mu = 0 \), the gradient \( \dot{v}(0, \sigma)=0 \), and so is the Hessian, except for the the first entry of \( \ddot{v}  \), that is,  %
\begin{align}
\tfrac{ \partial^{2}}{\partial \mu^{2}} v(\mu, \sigma) \Big |_{\mu = 0} = \frac{ \kappa_{1}^{2} - 2 \kappa_{0}^{2}}{\kappa_{0}^{2} \kappa_{1}^{2} \sigma^{2}} . 
\end{align}
As such, \( \log \BF_{1 \widetilde{0}} (Y^{n}) \) has a shifted asymptotically \( \chi^{2}(1) \)-distribution, i.e.,  
\begin{align}
n \Big ( \log \BF_{1 \widetilde{0}}(Y^{n} \, ; \, \kappa_{0}, \kappa_{1}) - \log \big ( \tfrac{\pi ( \theta \, | \, \Hc_{1})}{\pi( \theta \, | \, \Hc_{\widetilde{0}})} \big ) - E(\theta, n) \Big ) \overset{\Pf_{0, 1}}{\rightsquigarrow} \frac{ \kappa_{1}^{2} - 2 \kappa_{0}^{2}}{2 \kappa_{0}^{2} \kappa_{1}^{2}} Z^{2} ,
\end{align}
where \( Z \sim \Nc ( 0, 1) \). 

\begin{figure}[h]
\begin{tabular}{cc}
    \begin{minipage}{.5 \textwidth}
    \centering
    \qquad \scalebox{1.25}{\( \kappa_{0}=0.05 \)}
    \includegraphics[width= \linewidth]{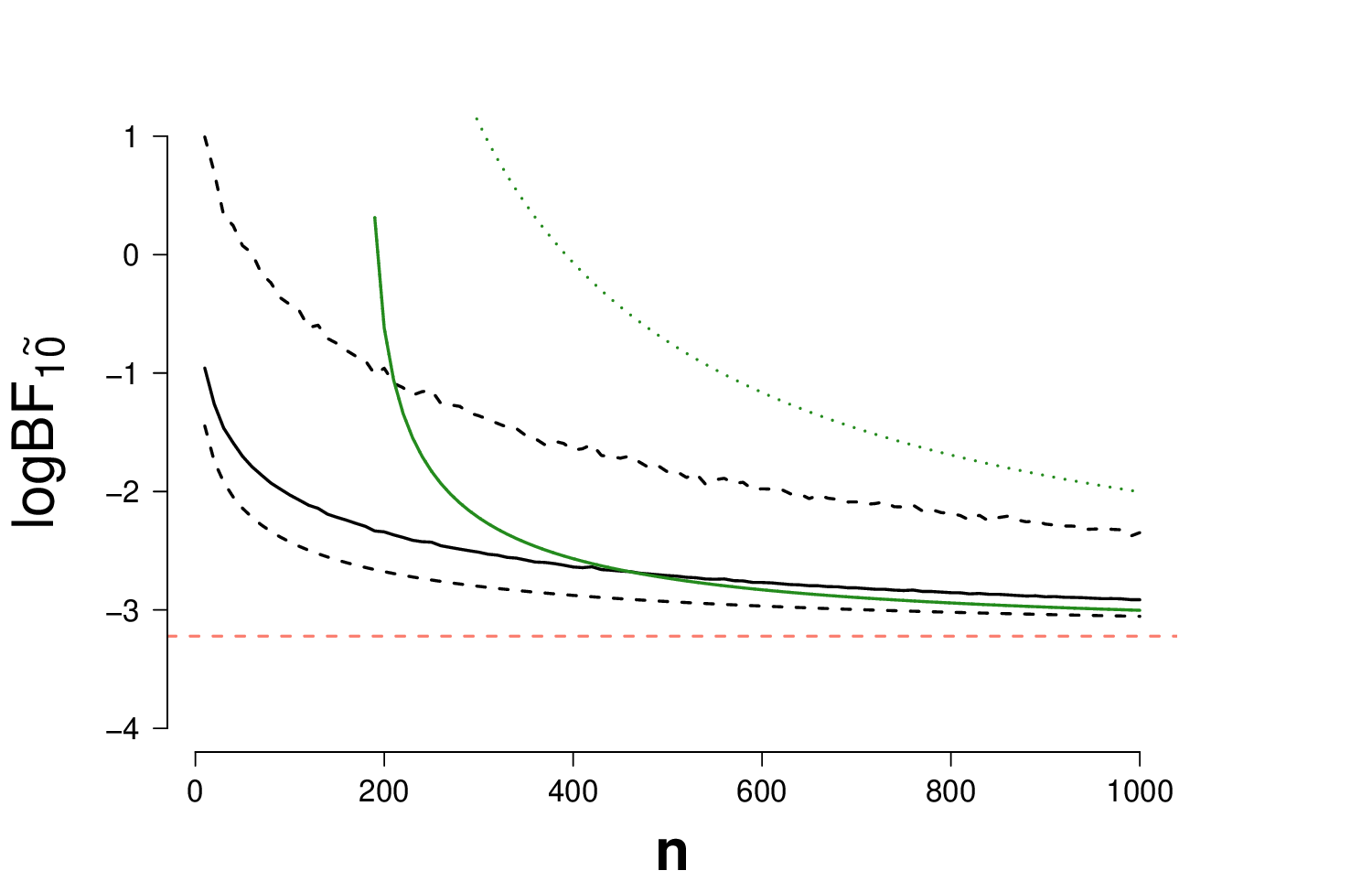}
    \end{minipage}  & 
    \begin{minipage}{.5 \textwidth}
    \centering
    \qquad \scalebox{1.25}{\( \kappa_{0}=0.10 \)}
    \includegraphics[width=\linewidth]{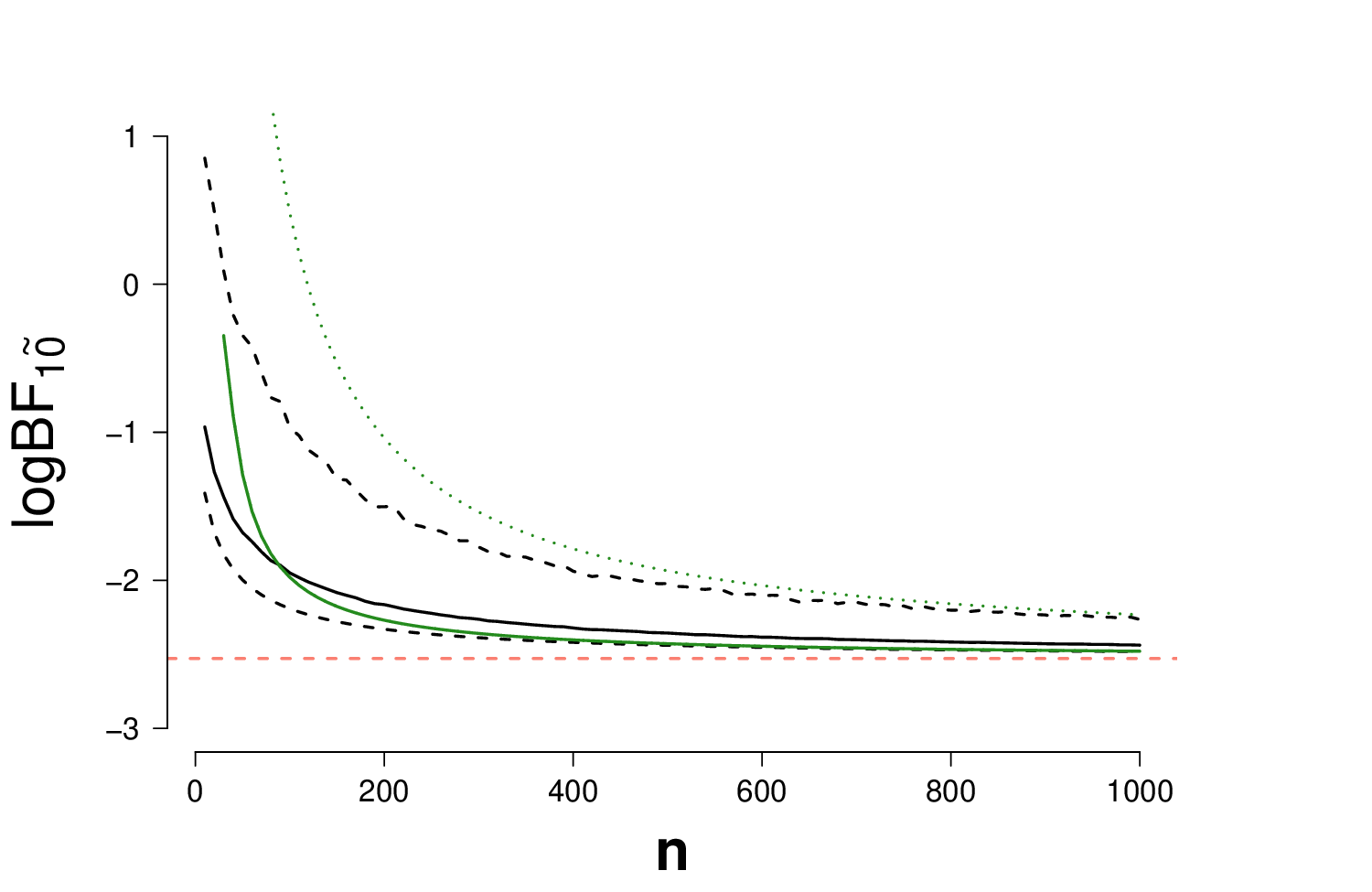}
    \end{minipage}  %
\end{tabular}
\caption{Under the null, the logarithm of the peri-null Bayes factor \( t \)-test has a shifted asymptotically \( \chi^{2}(1) \)-distribution with a mean (i.e., the solid curves) that decreases to the limit, e.g., \( \log \BF_{1 \widetilde{0}} = -3.22 \) and \( \log \BF_{1 \widetilde{0}} = -2.53 \) in the left and right plot respectively. The black and green curves correspond to the simulated and asymptotic \( \chi^{2}(1) \) sampling distribution respectively. The dotted curves show the 97.5\% and 2.5\% quantiles of the respective sampling distribution. Note that the convergence to the lower bound is slower when the peri-null is more concentrated, e.g., compare the left to the right plot.}
\label{figBfCNNull}
\end{figure}

More concretely, under \( \mu=0 \) and \( \sigma=1 \), \( \log  \BF_{1 \widetilde{0}}(Y^{n} \, ; \, 0.05, 1) \overset{\Pf_{0, 1}}{\rightarrow} -3.22 \), whereas \( \log  \BF_{1 \widetilde{0}}(Y^{n} \, ; \, 0.10, 1)  \) converges in probability to \( -2.53 \). Both cases yield evidence for the null hypothesis, but the evidence is stronger for the peri-null that is more tightly concentrated around 0. The approximation based on the asymptotic \( \chi^{2}(1) \)-distribution (in green) and the simulations (in black) are shown in \refFig{figBfCNNull}. 
In the left subplot, the curves based on the asymptotic \( \chi^{2}(1) \)-distribution only start from \( n=185 \), because only for \( n \geq 185  \) does \( \log ( 1 + C^{1}( 0, 1 \, | \, \Hc_{\widetilde{0}})/n + C^{2}( 0, 1 \, | \, \Hc_{\widetilde{0}})/n^{2}) \) have a non-negative argument; for \( \kappa_{0} = 0.05 \) we have that \(  C^{1}( 0, 1 \, | \, \Hc_{\widetilde{0}}) = -199.83 \). Note that under the null hypothesis, the Laplace approximations are accurate sooner than under the alternative hypothesis, because the priors are already concentrated at zero. Under the null hypothesis the general observation remains true that for reasonable sample sizes the expected peri-null Bayes factor is far from the limiting value. 

\begin{figure}[h]
\begin{tabular}{cc}
    \begin{minipage}{.5 \textwidth}
    \centering
    \qquad \scalebox{1.25}{\( \kappa_{0}=0.05 \)}
    \includegraphics[width= \linewidth]{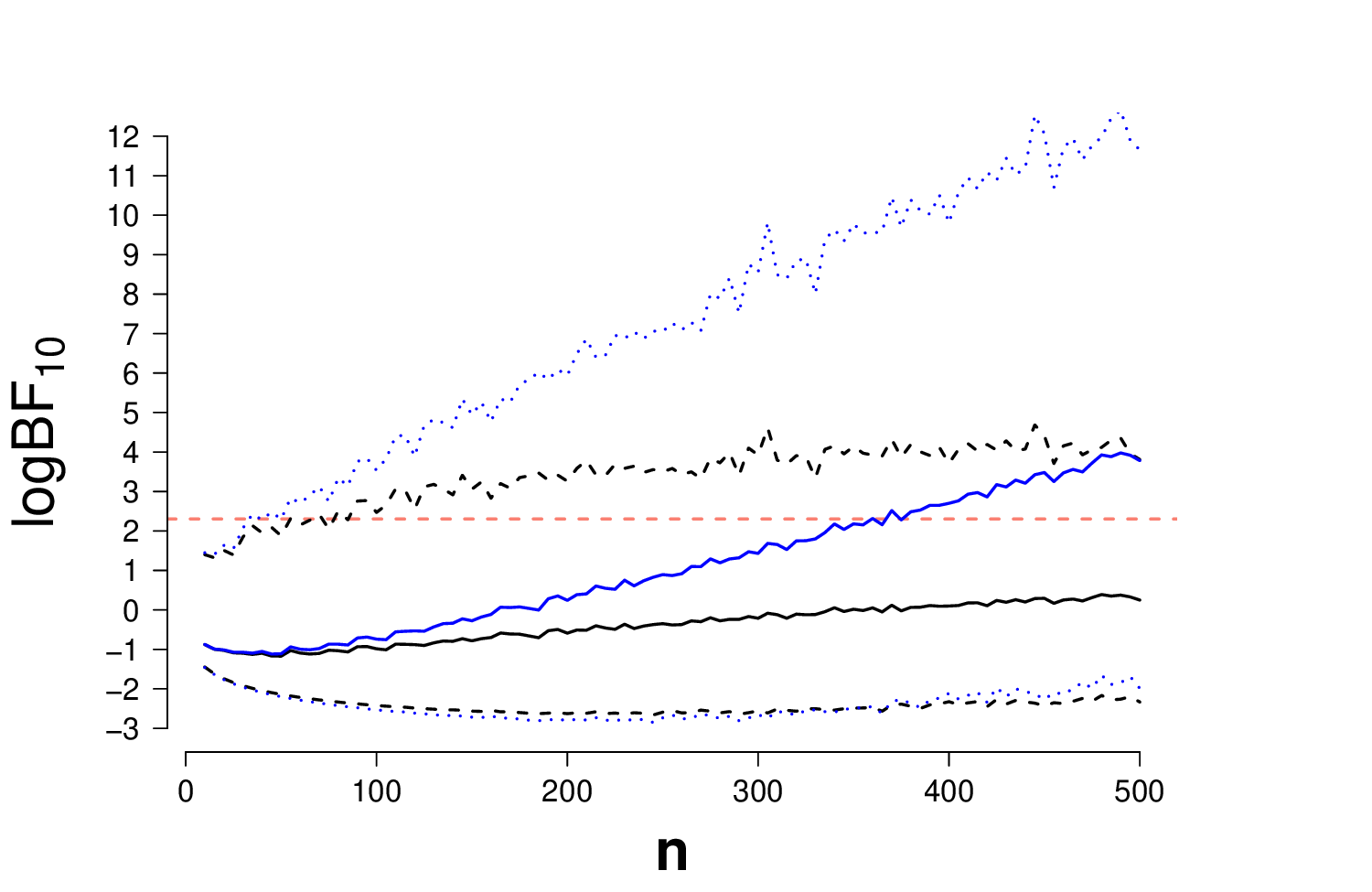}
    \end{minipage}  & 
    \begin{minipage}{.5 \textwidth}
    \centering
    \qquad \scalebox{1.25}{\( \kappa_{0}=0.10 \)}
    \includegraphics[width=\linewidth]{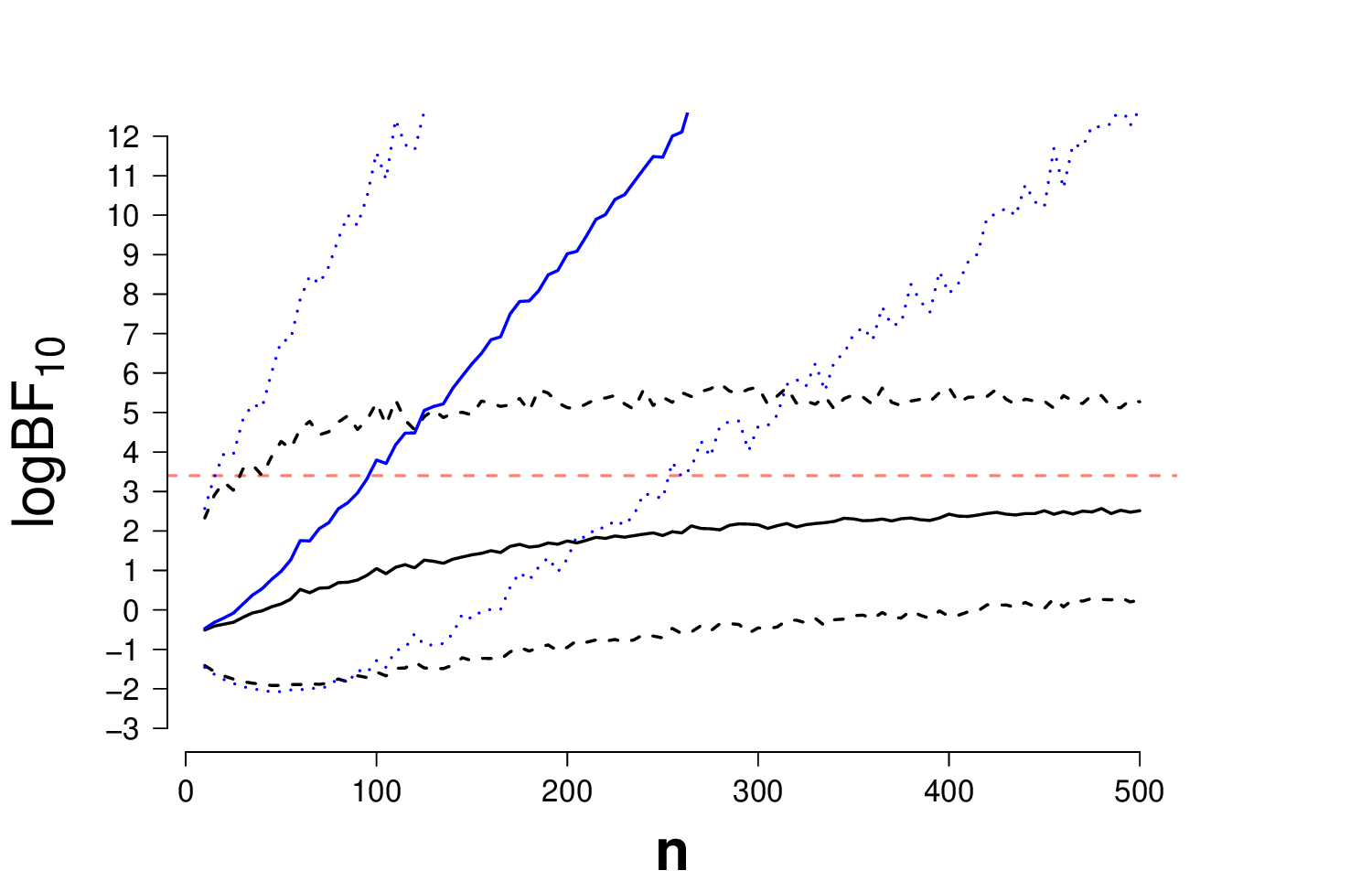}
    \end{minipage}  %
    \\
    \begin{minipage}{.5 \textwidth}
    \centering
    \includegraphics[width= \linewidth]{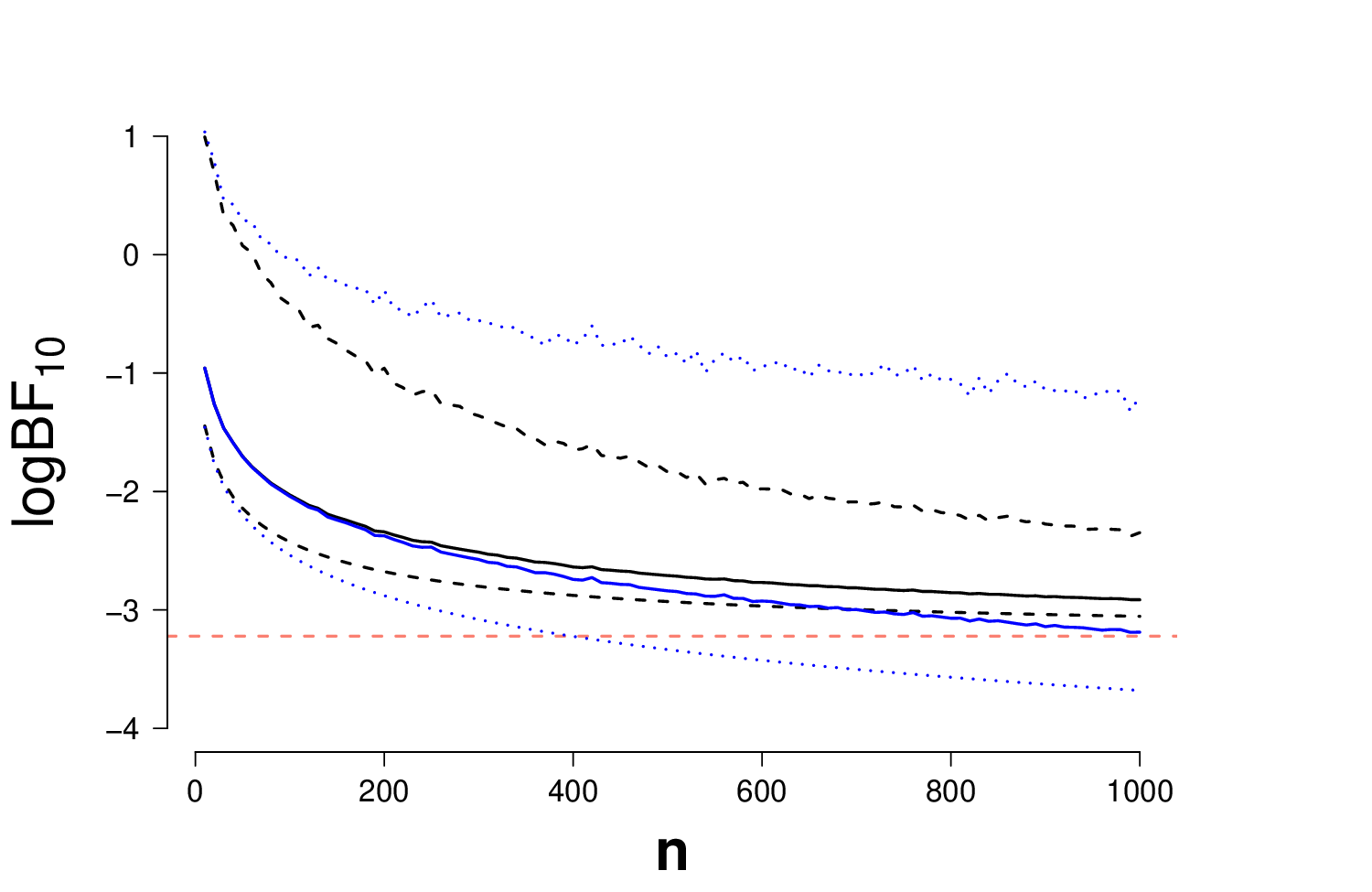}
    \end{minipage}  & 
    \begin{minipage}{.5 \textwidth}
    \centering
    \includegraphics[width=\linewidth]{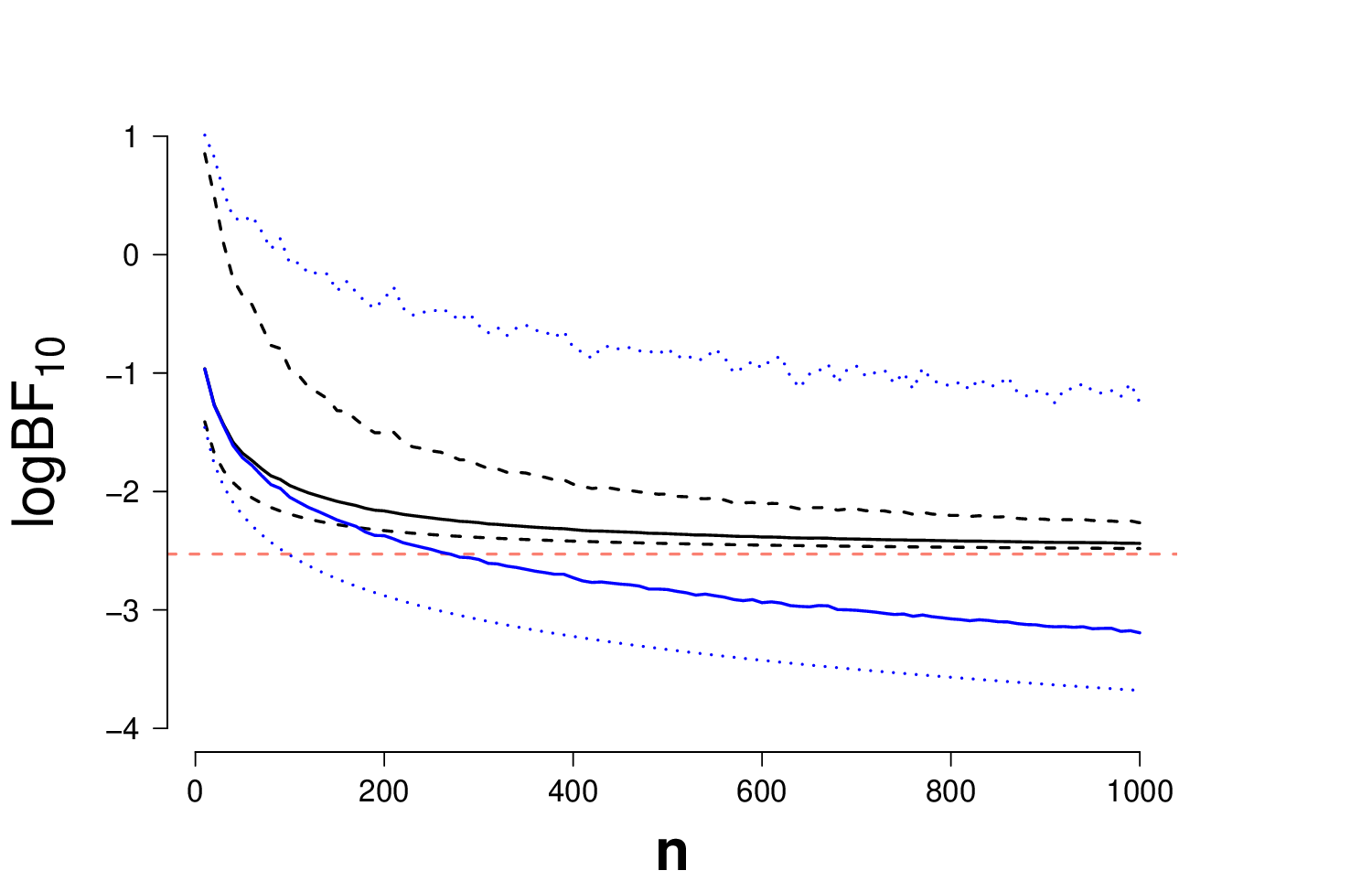}
    \end{minipage}  %
\end{tabular}
\caption{(Default) point-null Bayes factor \( t \)-tests (depicted in blue) are consistent under both the alternative and null, e.g., top and bottom row respectively, as opposed to peri-null Bayes factors (depicted in black). Note that the peri-null and the default point-null Bayes factors behave similarly when \( n \) is small. The domain where the two types of Bayes factors behave similarly is smaller when the peri-null is less concentrated, e.g., compare the right to the left column.}
\label{figBfCNCompareBfC0}
\end{figure}

Unlike the peri-null Bayes factor, the (default) point-null Bayes factor is consistent. \refFig{figBfCNCompareBfC0} shows the simulated sampling distribution of the point-null and peri-null Bayes factors in blue and black respectively. As before the 97.5\% quantile (top dotted curve), the average (solid curve), and the 2.5\% quantile (bottom dotted curve) are depicted as well. %

The top left subplot of \refFig{figBfCNCompareBfC0} shows that under \( \mu=0.167 \) and \( \sigma=1 \) the point-null and peri-null Bayes factor behave similarly up to \( n=30 \). Furthermore, the average point-null log Bayes factor crosses the peri-null upper bound of \( \log(10) \) at around \( n=380 \), whereas the peri-null Bayes factor remains bounded even in the limit, and is therefore inconsistent. The top right subplot shows, under \( \mu=0.348 \) and \( \sigma=1 \), that the discrepancy between the point-null and peri-null Bayes factor becomes apparent sooner when the peri-null prior is less concentrated, i.e., \( \kappa_{0} = 0.10 \) instead of \( \kappa_{0}=0.05 \). Also note that under these alternatives, the logarithm of the point-null Bayes factor grows linearly (e.g., \citealp{bahadur2009optimality, johnson2010use}). Hence, the point-null Bayes factor has a larger power to detect an effect than that afforded by the peri-null Bayes factor.

The bottom row of \refFig{figBfCNCompareBfC0} paints a similar picture; under the null the point-null Bayes factor accumulates evidence for the null hypothesis without bound as \( n \) grows. For \( \kappa_{0} = 0.05 \) the behavior of the peri-null and the point-null Bayes factor is similar up to \( n=200 \) and it takes about \( n=1,000 \) samples before the average point-null log Bayes factor crosses the peri-null lower bound of \( -3.22 \). For \( \kappa_{1} = 0.10 \) only \( n=270 \) samples are needed before the log Bayes factor for the point-null hypothesis crosses the peri-null lower bound of \( -2.53 \). 

\section{Towards Consistent Peri-Null Bayes Factors}
There are at least three methods to adjust the peri-null Bayes factor in order to avoid inconsistency. The first method changes both the point-null hypothesis \( \Hc_{0} \) and the alternative hypothesis \( \Hc_{1} \). Specifically, one may define the hypotheses under test to be non-overlapping (e.g., \citealp{BergerDelampady1987,ChandramouliShiffrin2019,rousseau2007approximating}). The resulting procedure is usually known as an `interval-null hypothesis', where the interval-null is defined as a (renormalized) slice of the prior distribution for the test-relevant parameter under an alternative hypothesis (e.g., \citealp{MoreyRouder2011}). For instance, in the case of a \( t \)-test an \emph{encompassing} hypothesis \( \Hc_{e} \) may assign effect size \( \delta \) a Cauchy distribution with location parameter 0 and scale \( \kappa_{e} \); from this encompassing hypothesis one may construct two rival hypotheses by restricting the Cauchy prior to particular intervals: the interval-null hypothesis truncates the encompassing Cauchy distribution to an interval centered on \( \delta=0 \): \( \delta \sim \text{Cauchy}(0, \kappa_{e})I(-a, a) \), whereas the interval-alternative hypothesis is the conjunction of the remaining two intervals, \( \delta \sim \text{Cauchy}(0, \kappa_{e})I(-\infty,-a) \) and \( \delta \sim \text{Cauchy}(0, \kappa_{e})I(a, \infty) \). As a consequence of \refThm{thmLimitPerinullBf}, or Theorem 1 (ii) of \citet{dawid2011posterior}, the resulting peri-null Bayes factor is then consistent in accordance to subjective interval belief; for all data-governing parameters \( \delta \) in the interior of the interval-null, \( \lim_{n \rightarrow \infty} \BF_{1 \widetilde{0}} = 0 \), and for \( \delta \) in the interior of the sliced out alternative \( \lim_{n \rightarrow \infty} \BF_{\widetilde{0}1} = 0  \).%
\footnote{For consistency to hold the standard condition is assumed that the interval-null or sliced up prior assigns positive mass to a neighborhood of \( \delta \) in the respective intervals.} %
In particular, when \( a = 1 \) and the data-governing \( \delta = 0.7 \), then this Bayes factor will eventually show unbounded evidence for the interval-null. %

One disadvantage of this method is the need to specify the width of the interval \citep[p. 367]{Jeffreys1961}. This disadvantage can be mitigated by reporting a range of non-overlapping interval-null Bayes factors as a function of \( a \); the researcher can then draw their own conclusion. The resulting range of interval-null Bayes factors also respects the uncertainty about the proper specification of the interval-null hypothesis and thereby avoids a false sense of precision.%
\footnote{We thank the first anonymous reviewer for suggesting this procedure to circumvent a definite choice for \( a \).} %
A second disadvantage of the non-overlapping interval-null method is that the prior distributions for the rival interval hypotheses are of an unusual shape -- a continuous distribution up to the point of truncation, where the prior mass abruptly drops to zero. It is debatable whether such artificial forms would ever result from an elicitation effort. A third disadvantage is that it seems somewhat circuitous to parry the critique ``the null hypothesis is never true exactly'' by adjusting both the null hypothesis \emph{and} the alternative hypothesis. 

The second method to specify a (partially) consistent peri-null Bayes factor is to change the point-null hypothesis to a peri-null hypothesis by supplementing rather than supplanting the spike with a distribution \citep{MoreyRouder2011}. In other words, the point-null hypothesis is upgraded to include a narrow distribution around the spike. This mixture distribution is generally known as a `spike-and-slab' prior, but here the slab represents the peri-null hypothesis and is relatively peaked. This mixture model \( \Hc_{0'} \) may be called a `hybrid null hypothesis' \citep{MoreyRouder2011}, a `mixture null hypothesis', or a `peri-point null hypothesis'. Thus, \( \Hc_{0'} = \xi \Hc_{0} + (1-\xi) \Hc_{\widetilde{0}} \), with \( \xi \in (0,1) \) the mixture weights and, say, \( \xi=\frac{1}{2} \). %
Because \( \xi > 0 \), the Bayes factor comparing \( \Hc_{0'} \) to \( \Hc_{1} \) will be consistent when the data come from \( \Hc_{0} \); and because \( \Hc_{0'} \) also has mass away from the point under test, the presence of a tiny true non-zero effect will not lead to the certain rejection of the null hypothesis as \( n \) grows large. The data determine which of the two peri-point components receives the most weight. %
As before, for modest sample sizes and small \( \kappa_{0} \), the distinction between point-null, peri-null, and peri-point null is immaterial. The main drawback of the peri-point null hypothesis is that it is consistent only when the data come from \( \Hc_{0} \); when the data come from \( \Hc_{1} \) or \( \Hc_{\widetilde{0}} \), the Bayes factor remains bounded as before (i.e., Eq.~\ref{eqPeriNullBf}).  

The third method is to define a peri-null hypothesis whose width \( \kappa_{0} \) slowly decreases with sample size (i.e., a `shrinking peri-null hypothesis'). For the \( t \)-test, one can take \( \kappa_{0} = c \sigma / \sqrt{n} \) for some constant \( c > 0 \) as proposed by \citet{BergerDelampady1987}, see also \citet{rousseau2007approximating}, except that their proposal involves a test of non-overlapping hypotheses. More generally, consistency follows by adapting \refThm{thmLimitPerinullBf}, which depends on the Laplace approximation that becomes invalid if \( \kappa_{0} \) shrinks too quickly. The representation \refEq{eq:BFcorrectionfactor} shows that this consistency fix is equivalent to keeping the peri-null correction Bayes factor \( \BF_{\widetilde{0} 0} \) close to one regardless of the data. Note that this is attainable as \( \kappa_{0} \rightarrow 0 \) the peri-null and point-null become indistinguishable. In other words, the consistency of such a shrinking peri-null Bayes factor is essentially driven by the asymptotic behavior of the point-null Bayes factor; arguably, one might as well employ this point-null Bayes factor to begin with. One drawback of the shrinking peri-null hypothesis is that it is incoherent, because the prior distribution depends on the intended sample size.%
\footnote{We term a Bayes factor incoherent if the result depends on whether the data are analyzed all at once or batch by batch.} %
There could nevertheless be a pragmatic argument for tailoring the definition of the peri-null to the resolving power of an experiment.

\section{Concluding Comments}
The objection that ``the null hypothesis is never true'' may be countered by abandoning the point-null hypothesis in favor of a peri-null hypothesis. For moderate sample sizes and relatively narrow peri-nulls, this change leaves the Bayes factor relatively unaffected. For large sample sizes, however, the change exerts a profound influence and causes the Bayes factor to be inconsistent, with a limiting value given by the ratio of prior ordinates evaluated at the maximum likelihood estimate (cf. \citealp[p. 367]{Jeffreys1961} and \citealp[pp. 411-412]{MoreyRouder2011}). Hence, we believe that as far as Bayes factors are concerned, there is much to lose and little to gain from adopting a peri-null hypothesis in lieu of a point-null hypothesis. %
Here we also derived the asymptotic sampling distribution of the peri-null Bayes factor and show that its limiting value is essentially an upper bound under the alternative and a lower bound under the null. The asymptotic distributions also provide insights to typical values of the peri-null Bayes factor at a finite \( n \).  
Inconsistency may not trouble subjective Bayesians: if the peri-null hypothesis truly reflects the belief of a subjective sceptic, and the alternative hypothesis truly reflects the belief of a subjective proponent, then the Bayes factor provides the relative predictive success for the sceptic versus the proponent, and it is irrelevant whether or not this relative success is bounded. Objective Bayesians, however, develop and apply procedures that meet various desiderata (e.g., \citealp{BayarriEtAl2012,ConsonniEtAl2018}), with consistency a prominent example. As indicated above, the desire for consistency was the primary motivation for the development of the Bayesian hypothesis test \citep{WrinchJeffreys1921}. For objective Bayesians then, it appears the point-null hypothesis is more than just a mathematically convenient approximation to the peri-null hypothesis \citep[p. 367]{Jeffreys1961}. The peri-point mixture model (consistent only under the point-null hypothesis) and the shrinking peri-point model (incoherent because the prior width depends on sample size) may provide acceptable compromise solutions.

Regardless of one's opinion on the importance of consistency, it is evident that seemingly inconsequential changes in prior specification may asymptotically yield fundamentally different results. Researchers who entertain the use of peri-null hypotheses should be aware of the asymptotic consequences; in addition, it generally appears prudent to apply several tests and establish that the conclusions are relatively robust. 

\section*{Acknowledgements}
This research was supported by the Netherlands Organisation for Scientific Research (NWO; grant \#016.Vici.170.083). The authors would like to thank Paulo Serra and Muriel P\'{e}rez for their enlightening discussions regarding the asymptotics discussed in this paper. We are also grateful for the constructive suggestions for improvement by the editor and two anonymous reviewers.

\bibliographystyle{abbrvnat}
\bibliography{Alexander,referenties}

\newpage

\appendix
\section{Laplace Approximation}
\label{appLaplace}

The Laplace approximation uses a (multivariate) Taylor expansion for which we introduce notation. Let \( h : \Theta \subset \R^{p} \rightarrow \R \), i.e., \( h(\theta)= - \tfrac{1}{n} \sum_{i=1}^{n} \log f(y_{i} \, | \, \theta ) \), and we write \( \hat{\theta} \) for the point in its domain where \( h \) takes its global minimum. Furthermore, we use subscripts to denote partial derivatives, whereas superscripts refer to components of a vector, or more generally an array. For instance, \( \pi_{a} = \tfrac{ \partial}{\partial \theta_{a}} \pi(\hat{\theta})  \) refers to the \( a \)-th component of the vector of partial derivatives \( [D^{1} \pi(\hat{\theta})] \) of the prior \( \pi \) evaluated at the MLE. Similarly, we write \( h_{abc} = \tfrac{\partial^{3}}{ \partial \theta_{a} \partial \theta_{b} \partial \theta_{c} } h(\hat{\theta}) \) for the \( abc \)-th component of the three-dimensional array \( [D^{3} h(\hat{\theta}) ] \). Hence, the number of indices in the subscript corresponds to the number of derivatives of \( h \) and the indices, each in \( 1, 2, \ldots, p \), provide the location of the component. 

We use superscripts to refer to the component of a vector. For instance, \( \tilde{q}^{a} = (\theta^{a} - \hat{\theta}^{a}) \) represents the \( a \)-th component of the difference vector \( \tilde{q} = \theta - \hat{\theta} \), thus, equivalently \( \tilde{q}^{a} := e_{a}^{T} \tilde{q} \), where \( e_{a} \) is the unit (column) vector with entry \( 1 \) at index \( a \) and zero elsewhere. Similarly, \( \varsigma^{abcd} \) the \( abcd \)-th component of a four dimensional array. 

Moreover, we employ Einstein's summation convention and suppress the sum whenever an index occurs in both the sub and superscript. For instance, 
\begin{align}
h_{a} \tilde{q}^{a} &:= \sum_{a=1}^{p} h_{a} \tilde{q}^{a}, \\
h_{abc} \tilde{q}^{a} \tilde{q}^{b} \tilde{q}^{c} & := \sum_{a=1}^{p} \sum_{b=1}^{p} \sum_{c=1}^{p} h_{abc} \tilde{q}^{a} \tilde{q}^{b} \tilde{q}^{c}.
\end{align}
The former defines an inner product between the gradient of \( h \) and deviations \( \tilde{q} \), whereas the \( h_{abc} = [D^{3} h]_{abc} \) refers to the \( a \)-th row, \( b \)-th column, and \( c \)-th depth of the three-dimensional array consisting of partial derivatives of \( h \) of order three. Lastly, we use the shorthand notation 
\begin{align}
h_{a} h_{b} \tilde{q}^{a} \tilde{q}^{b} := \sum_{a} \sum_{b} h_{a} h_{b} \tilde{q}^{a} \tilde{q}^{b},
\end{align}
to denote the nested sum which is needed for Cauchy products \( (h_{a} \tilde{q}^{a}) (h_{b} \tilde{q}^{b}) \). For instance, with \( d =2 \)
\begin{align}
(h_{1} \tilde{q}^{1} + h_{2} \tilde{q}^{2}) (h_{1} \tilde{q}^{1} + h_{2} \tilde{q}^{2}) = h_{1} \tilde{q}^{1} h_{1} \tilde{q}^{1} + 2 h_{1} \tilde{q}^{1} h_{2} \tilde{q}^{2} + h_{2} \tilde{q}^{2} h_{2} \tilde{q}^{2} ,
\end{align}
which is equivalent to 
\begin{align}
h_{1} h_{1} \tilde{q}^{1} \tilde{q}^{1} + h_{1} h_{2} \tilde{q}^{1} \tilde{q}^{2} +h_{2} h_{1} \tilde{q}^{2} \tilde{q}^{1} + h_{2} h_{2} \tilde{q}^{2} \tilde{q}^{2} .
\end{align}
With these notational conventions a multivariate Taylor approximation is denoted as %
\begin{align}
h(\theta)  & = h ( \hat{\theta}) + h_{a} \tilde{q}^{a} + \tfrac{h_{ab}}{2!}  \tilde{q}^{a} u^{b} + \tfrac{h_{abc}}{3!}   \tilde{q}^{a} \tilde{q}^{b} \tilde{q}^{c} + \tfrac{h_{abcd}}{4!}   \tilde{q}^{a} \tilde{q}^{b} \tilde{q}^{c} \tilde{q}^{d}  + \Oc( | u |^{5} ). %
\end{align}
and note the similarity to its one-dimensional counterpart. 

\begin{theorem}[Laplace expansion with error term]
Let \( \Pc_{\Theta} \) be a collection of density functions that are six times continuously differentiable in \( \theta \in \Theta \subset \R^{p} \), and \( \pi(\theta) \) a prior density that is four times continuously differentiable. Let \( Y \iidSim f(y \, | \, \theta) \) for certain \( \theta \), then with  \( \hat{\theta} \) the MLE %
\begin{align}
p(y^{n}) & = \int_{\Theta} f(y^{n} \, | \, \theta) \pi(\theta) \der \theta \\
& = (\tfrac{2 \pi}{n})^{\tfrac{p}{2}} f(y^{n} \, | \, \hat{\theta}) \pi ( \hat{\theta}) | \hat{I} ( \hat{\theta})|^{-1/2} \Big [1 + \tfrac{C^{(1)}(\hat{\theta})}{n} + \tfrac{C^{(2)}(\hat{\theta})}{n^{2}} + \Oc ( n^{-3}) \Big ] ,
\end{align}
where \( | \cdot | \) denotes the determinant and 
\begin{align}
\label{laplaceOrder1}
C^{(1)}(\hat{\theta}) & = \tfrac{\pi_{ab}}{2 \pi(\hat{\theta})}  \varsigma^{ab} - \big ( \tfrac{h_{abcd} }{24} +\tfrac{h_{abc}  \pi_{u}  }{6 \pi(\hat{\theta})} \big ) \varsigma^{abcd} + \tfrac{h_{abc} h_{uef} }{72}  \varsigma^{abcdef} , \\
\label{laplaceOrder2}
C^{(2)}(\hat{\theta}) & = \tfrac{ \pi_{abcd}}{24 \pi(\hat{\theta})} \varsigma^{abcd} - \tfrac{\pi(\hat{\theta}) h_{abcdef} + 6 h_{abcde} \pi_{f} +15 h_{abcd} \pi_{ef} +20 h_{abc} \pi_{def} }{720 \pi(\hat{\theta})} \varsigma^{abcdef} \\
\nonumber
& + \tfrac{ 5 \pi(\hat{\theta}) h_{abcd} h_{efgh} + 8 \pi(\hat{\theta})  h_{abcde} h_{fgh} + 40 h_{abc} \big ( h_{defg} \pi_{h} + h_{def} \pi_{gh} \big) }{5760 \pi(\hat{\theta})} \varsigma^{abcdefgh} \\
\nonumber
& - \tfrac{3 \pi(\hat{\theta}) h_{abcd} h_{efg} h_{hij} + 4 h_{abc} h_{def} h_{ghi} \pi_{j}}{5184 \pi(\hat{\theta})} \varsigma^{abcdefghij} + \tfrac{ h_{abc} h_{def} h_{ghi} h_{jkl}}{31104} \varsigma^{abcdefghijkl}
\end{align}
where \( \varsigma^{ab}  \), \( \varsigma^{abcd} \), \( \varsigma^{abcdef} \), \( \varsigma^{abcdefgh} \),  \( \varsigma^{abcdefghij} \), and \( \varsigma^{abcdefghijkl} \) represent the \( ab \)-th component of the second, the \( abcd \)-th component of the fourth, the \( abcdef \)-th component of the sixth moment, the \( abcdefgh \)-th component of the eigth moment, the \( abcdefghij \)-th component of the tenth moment, and the \( abcdefghijkl \)-th component of the twelfth moment, of the \( p \) dimensional random vector \( Q \sim \Nc_{p} ( 0, \hat{I}(\hat{\theta})^{-1}) \), respectively. \( \hfill \diamond \)
\end{theorem}

\begin{proof}
The proof is based on (i) Taylor-expanding the exponential of the log-likelihood of order five around \( \hat{\theta} \), (ii) the definition of the exponential as a series and Taylor-expanding \( \pi \) to third order at the same point \( \hat{\theta} \), and (iii) properties of the normal distribution. 

\textbf{Step (i)} Let \( h(\theta) = \tfrac{1}{n} \sum_{i=1}^{n} \log f(y_{i} \, | \, \theta) \), then since \( h(\theta) \in C^{6}(\Theta) \) we know that there exists \( \delta > 0  \) such that in a ball \( B_{\hat{\theta}}(\delta) \subset \R^{p} \) of radius \( \delta \) centered at \( \hat{\theta} \) the average log-likelihood \( h_{n}(\theta) \) is well-approximated by a Taylor expansion of order 5. This combined with \( \hat{\theta} \) being the MLE and the notation \( \tilde{q} = \theta - \hat{\theta} \) yields %
\begin{align}
p(y^{n}) = & \int_{\Theta} e^{ - n h (\theta)} \pi(\theta) \der \theta = \int_{ B_{\hat{\theta}}(\delta)} e^{ - n h(\hat{\theta}) - \tfrac{n
h_{ab} \tilde{q}^{a} \tilde{q}^{b}}{2} + \tilde{R}(\tilde{q})} \pi(\tilde{q}) \der \tilde{q} , \\
= & f(y^{n} \, | \, \hat{\theta})  \int_{ B_{\hat{\theta}}(\delta) }   e^{- \tfrac{n
h_{ab} \tilde{q}^{a} \tilde{q}^{b}}{2!}} e^{- \tilde{R}(\tilde{q})} \pi(\tilde{q}) \der \tilde{q} ,
\end{align}
where %
\begin{align}
\tilde{R}(\tilde{q}) =  n [\tfrac{h_{abc} \tilde{q}^{a} \tilde{q}^{b} \tilde{q}^{c}}{3!} + \tfrac{h_{abcd} \tilde{q}^{a} \tilde{q}^{b} \tilde{q}^{c} \tilde{q}^{d}}{4!} + \tfrac{h_{abcde} \tilde{q}^{a} \tilde{q}^{b} \tilde{q}^{c} \tilde{q}^{d} \tilde{q}^{e}}{5!} + \Oc( |\tilde{q}|^{6}) ]  ,
\end{align}
is the bounded remainder term since \( h \in C^{6}(\Theta) \). The replacement of \( \Theta \) by \( B_{\hat{\theta}}(\delta) \) in the integral is justified if the mass is concentrated at \( \hat{\theta} \), thus, whenever the integral with respect to the first order term falls off quadratically, that is, if 
\begin{align}
|n \hat{I}(\hat{\theta})|^{1/2} e^{-n (h(\theta) - h(\hat{\theta}))} \pi(\theta) 
= \Oc(n^{-2}) , 
\end{align}
which is the case when \( \hat{\theta} \) is unimodal. When it is not unimodal, but \( \hat{\theta} \) is a global maximum, then the condition implies that the requirement that the contribution of the other maxima is not too big.

\textbf{Step (ii)} After centering the integral at \( \hat{\theta} \) we scale with respect to \( \sqrt{n} \), that is, we apply the change of variable \( q =  \sqrt{n} \tilde{q} \), thus, \( \int n^{-p/2} \der q = \int \der \tilde{q} \) and therefore
\begin{align}
p(y^{n}) & = (\tfrac{2 \pi}{n})^{p/2} f(y^{n} \, | \, \hat{\theta}) | \hat{I} ( \hat{\theta}) |^{-1/2} \int_{ B_{\hat{\theta}}(\sqrt{n} \delta) } \tilde{\varphi} ( q) e^{-R(q)} \tilde{\pi}(q) \der q ,
\end{align}
where \( \tilde{\varphi} \) is the density of a multivariate normal distribution centered at \( 0 \) and covariance matrix \( \Sigma = \hat{I}^{-1}(\hat{\theta}) \), and where \( \tilde{\pi}(q) \) is the Taylor approximation of \( \pi \) at the MLE, that is, 
\begin{align}
\tilde{\pi}(q) = \pi(\hat{\theta}) + \tfrac{\pi_{a} (\hat{\theta}) q^{a} }{ n^{1/2} } + \tfrac{ \pi_{ab}(\hat{\theta}) q^{a} q^{b} }{2! n} + \tfrac{\pi_{abc} (\hat{\theta}) q^{a} q^{b} q^{c} }{3! n^{3/2}} + \Oc(n^{-2}) , 
\end{align}
and where the remainder term is now %
\begin{align}
\nonumber
R(q) =  \tfrac{h_{abc} q^{a} q^{b} q^{c}}{3! n^{1/2} } + \tfrac{h_{abcd} q^{a} q^{b} q^{c} q^{d}}{4! n} + \tfrac{h_{abcde} q^{a} q^{b} q^{c} q^{d} q^{e}}{5! n^{3/2} } + \tfrac{h_{abcde} q^{a} q^{b} q^{c} q^{d} q^{e} q^{f} }{6! n^{2} } \Oc( n^{-5/2})  . 
\end{align}
To exploit the properties of Gaussian integrals we replace integration domain \( B_{\hat{\theta}}(\sqrt{n} \delta) \) by \( \R^{p} \), which is justified when \( n \) is large, and because the tails of a normal density fall off exponentially. 

By definition of \( e^{-R(q)} \) as a series and without the exponential approximation error %
\begin{align}
\label{eqLaplaceExpansion1}
p(y^{n}) & \approx (\tfrac{2 \pi}{n})^{p/2} f(y^{n} \, | \, \hat{\theta}) | \hat{I} ( \hat{\theta}) |^{-1/2} \\
\nonumber
& \times \int_{ \R^{p} } \tilde{\varphi} ( q) \big [1 - R(q) + \tfrac{ R(q)^{2}}{2!} - \tfrac{R(q)^{3}}{3!} + \Oc ( |R(q)|^{4}) \big ] \tilde{\pi}(q) \der q .
\end{align}
From here onwards we focus on the integral \refEq{eqLaplaceExpansion1}, which after some straightforward but tedious computations can be shown to be of the form %
\begin{align}
\int_{\R^{p}} \tilde{\varphi}( q ) \big [A_{0} + A_{1} n^{-1/2} + A_{2} n^{-1} + A_{3} n^{-3/2} + A_{4} n^{-2} +\Oc( n^{-3}) \big ] \der q ,
\end{align}
where the \( A_{j} \) terms are functions of \( q \) and \( \hat{\theta} \) defined by the series representation of \( e^{-R(q)} \) and \( \tilde{\pi}(q) \). 

\textbf{Step (iii)} The terms \( A_{j} \) are given below. Of the following results only the exact values of \( A_{0}, A_{2} \) and \( A_{4} \) matter; what matters for \( A_{1} \) and \( A_{3} \) is that they only involve odd powers of \( q \): %
\begin{align}
A_{0} & = \pi(\hat{\theta}) \\
A_{1} & =  \pi_{a}  q^{a} - \tfrac{h_{abc} \pi(\hat{\theta})}{6} q^{a} q^{b} q^{c}  \\
A_{2} & = \tfrac{\pi_{ab} }{2}  q^{a} q^{b} - \big ( \tfrac{\pi(\hat{\theta}) h_{abcd} }{24} +\tfrac{h_{abc} \pi_{u} }{6} \big ) q^{a} q^{b} q^{c} q^{d} + \tfrac{\pi(\hat{\theta}) h_{abc} h_{uef} }{72}  q^{a} q^{b} q^{c} q^{d} q^{e} q^{w} \\ 
A_{3} & = \tfrac{ \pi_{abc}}{6} q^{a} q^{b} q^{c} %
- \tfrac{6 h_{abcde} \pi(\hat{\theta}) + 30  h_{abcd} \pi_{v} + 60 h_{abc} \pi_{uv}}{720} q^{a} q^{b} q^{c} q^{d} q^{e}   \\
\nonumber
& + \tfrac{h_{abc} h_{uefl} \pi(\hat{\theta}) + 2 h_{abc} h_{uef} \pi_{l} }{144} q^{a} q^{b} q^{c} q^{d} q^{e} q^{w} q^{l} \\ 
A_{4} & = \tfrac{ \pi_{abcd}}{24} q^{a} q^{b} q^{c} q^{d} \\
\nonumber
& - \tfrac{\pi(\hat{\theta}) h_{abcdef} + 6 h_{abcde} \pi_{f} +15 h_{abcd} \pi_{ef} +20 h_{abc} \pi_{def} }{720 } q^{a} q^{b} q^{c} q^{d} q^{e} q^{f} \\
\nonumber
& + \tfrac{ 5 \pi(\hat{\theta}) h_{abcd} h_{efgh} + 8 \pi(\hat{\theta})  h_{abcde} h_{fgh} + 40 h_{abc} \big ( h_{defg} \pi_{h} + h_{def} \pi_{gh} \big) }{5760} q^{a} q^{b} q^{c} q^{d} q^{e} q^{f} q^{g} q^{h} \\
\nonumber
& - \tfrac{3 \pi(\hat{\theta}) h_{abcd} h_{efg} h_{hij} + 4 h_{abc} h_{def} h_{ghi} \pi_{j}}{5184} q^{a} q^{b} q^{c} q^{d} q^{e} q^{f} q^{g} q^{h} q^{i} q^{j} \\
\nonumber
& + \tfrac{\pi(\hat{\theta}) h_{abc} h_{def} h_{ghi} h_{jkl}}{31104} q^{a} q^{b} q^{c} q^{d} q^{e} q^{f} q^{g} q^{h} q^{i} q^{j} q^{k} q^{l} .
\end{align}
Since for \( k \) odd \( A_{k} \) only involve odd powers of \( q \) we conclude that their integral with respect to \( \tilde{\varphi}(q) \) vanishes. Hence, %
\begin{align}
p(y^{n}) =  (\tfrac{2 \pi}{n})^{p/2} f(y^{n} \, | \, \hat{\theta}) | \hat{I} ( \hat{\theta}) |^{-1/2} \pi(\hat{\theta}) \big [1 + \tfrac{E[A_{2}]}{n \pi(\hat{\theta})} +  \tfrac{E[A_{4}]}{n^{2} \pi(\hat{\theta})}  \Oc( n^{-3}) \big ]  ,
\end{align}
where \( E[A_{2}] \) and \( E[A_{4}] \) are expectations with respect to \( Q \sim \Nc ( 0, \hat{I}(\hat{\theta})^{-1}) \). This implies that the order \( n^{-1} \) and \( n^{-2} \) terms in the assertion are \( C^{(1)}(\hat{\theta}) = E[A_{2}]/\pi(\hat{\theta}) \) and \( C^{(2)}(\hat{\theta}) = E[A_{4}]/\pi(\hat{\theta}) \). %
\end{proof}

The components of higher moments can be expressed in terms of the covariances \( \varsigma^{ab} = \cov ( Q^{a}, Q^{b}) \) using Isserlis' formula (\citealp{isserlis1918formula, mccullagh2018tensor}). For moments \( \varsigma^{a_{1} \cdots a_{w}} \), that is, a component of the \( w \)th moment of \( Q \) with \( w = 2 v  \) even, the following holds
\begin{align}
\varsigma^{a_{1} \cdots a_{w} } = \sum_{u \in P_{w}^{2}} \prod_{i,j \in u} \varsigma^{ij} , 
\end{align}
where \(  P_{w}^{2} \) is the collection of all pairs of which there are \( v \). For instance, for \( w=4 \), \( \varsigma^{abcd} \) is a sum of 2-products of pairs, for \( w=6 \) is a sum of 3-products of \( \varsigma^{abcdef} \) and so forth and so on. More specifically, %
\begin{align}
\varsigma^{abcd} & = \varsigma^{ab} \varsigma^{cd}  + \varsigma^{ac} \varsigma^{bd}  + \varsigma^{ad} \varsigma^{bc} \\
\varsigma^{abcdef} & = \varsigma^{ab} \varsigma^{cd}  \varsigma^{ef} + \varsigma^{ab} \varsigma^{ce}  \varsigma^{df} + \varsigma^{ab} \varsigma^{cf}  \varsigma^{de} \\
\nonumber
& + \varsigma^{ac} \varsigma^{bd}  \varsigma^{ef} + \varsigma^{ac} \varsigma^{be}  \varsigma^{df} + \varsigma^{ac} \varsigma^{bf}  \varsigma^{de} \\
\nonumber
& + \varsigma^{ad} \varsigma^{bc}  \varsigma^{ef} + \varsigma^{ad} \varsigma^{be}  \varsigma^{cf} + \varsigma^{ad} \varsigma^{bf}  \varsigma^{ce} \\
\nonumber
& + \varsigma^{ae} \varsigma^{bc}  \varsigma^{df} + \varsigma^{ae} \varsigma^{bd}  \varsigma^{cf} + \varsigma^{ae} \varsigma^{bf}  \varsigma^{cd} \\
\nonumber
& + \varsigma^{af} \varsigma^{bc}  \varsigma^{de} + \varsigma^{af} \varsigma^{bd}  \varsigma^{ce} + \varsigma^{af} \varsigma^{be}  \varsigma^{cd} ,
\end{align}
where all indexes \( a,b,c,d,e,f =1,2, \ldots, p \). The expression of \( \varsigma^{abcdefgh} \), \( \varsigma^{abcdefghij} \), and \( \varsigma^{abcdefghijkl} \) define sums of \( 105=3 \times 5 \times 7 \), \( 945=3 \times 5 \times 7 \times 9 \), and \( 10,395=3 \times 5 \times 7 \times 9 \times 11 \) terms respectively, and due to space restrictions their exact forms are not displayed here. 

\end{document}